\documentclass[a4paper]{article}

\usepackage{amssymb,amsmath}
\usepackage[plainpages=false, hidelinks,
bookmarks]{hyperref}
\usepackage{graphicx}
\newtheorem{theo}{Theorem}
\newtheorem{lemma}{Lemma}
\newtheorem{rema}[theo]{Remark}
\begin{document}

\begin{center}
	{\Large \bf
		Tridiagonal test matrices for eigenvalue computations: two-parameter extensions of the Clement matrix\footnote{This is a preprint of a paper whose final and definite form is in Journal of Computational and Applied Mathematics: 
			
		R.~Oste, J.~Van~der~Jeugt, Tridiagonal test matrices for eigenvalue computations: Two-parameter extensions of the Clement matrix, Journal of Computational and Applied Mathematics, \textbf{314} (2017), 30--39, ISSN 0377-0427, doi: \href*{http://dx.doi.org/10.1016/j.cam.2016.10.019}.}
	} \\[5mm]
	{\bf Roy Oste, Joris Van der Jeugt}\\[3mm]	
	E-mail addresses: {\tt Roy.Oste@UGent.be},	{\tt Joris.VanderJeugt@UGent.be}\\[3mm]
	Department of Applied Mathematics, Computer Science and Statistics, Faculty of Sciences, Ghent University, Krijgslaan 281-S9, 9000 Gent, Belgium.
\end{center}

\begin{abstract}
	The Clement or Sylvester-Kac matrix is a tridiagonal matrix with zero diagonal and simple integer entries. Its spectrum is known explicitly and consists of integers which makes it a useful test matrix for numerical eigenvalue computations. 
	We consider a new class of appealing two-parameter extensions of this matrix which have the same simple structure and whose eigenvalues are also given explicitly by a simple closed form expression. 
	The aim of this paper is to present in an accessible form these new matrices and examine some numerical results regarding the use of these extensions as test matrices for numerical eigenvalue computations. 
\end{abstract}


	\noindent \emph{Keywords:}
	Tridiagonal matrix; Exact eigenvalues; Eigenvalue computation; Clement matrix; Sylvester-Kac matrix
	
\noindent 	\emph{2010 MSC:} 65F15; 15A18; 15A36
	
	
\section{Introduction}

For a positive integer $n$, consider the $(n+1)\times(n+1)$ matrix $C_n$ whose non-zero entries are given by
\begin{equation}
	\label{entries}
	c_{k,k+1} = c_{n+2-k,n+1-k} = k\qquad\mbox{for }k\in \{1,\dots,n\},
\end{equation}
or explicitly, the matrix 
\begin{equation}
C_n= 
\left(\begin{array}{ccccccc}
0 & 1 &  &   & &\\
n & 0 &  2& & &\\
& n-1& 0  & 3 & &\\
&   & \ddots &\ddots &\ddots  &  \\
&   &  &   3  & 0 & n-1&  \\
&   &  &     & 2 & 0& n \\
&   &  &  & & 1& 0
\end{array}\right).
\label{KSC}
\end{equation}
This matrix appears in the literature under several names: the Sylvester-Kac matrix, the Kac matrix, the Clement matrix, etc.
It was already considered by Sylvester~\cite{Sylvester}, used by M.~Kac in some of his seminal work~\cite{Kac}, proposed by Clement as a test matrix for eigenvalue computations~\cite{Clement}, and it continues to attract attention~\cite{Taussky, Boros, Bevilacqua, Edelman}. 

The matrix $C_{n}$
has a simple structure, it is tridiagonal with zero diagonal and has integer entries. The main property of $C_{n}$ is that its spectrum is known explicitly and is remarkably simple; the eigenvalues of $C_{n}$ are the integers
\begin{equation}
	-n,-n+2,-n+4,\dotsc,n-2,n.
	\label{KacEig}
\end{equation}
The $n+1$ distinct eigenvalues are symmetric around zero, equidistant and range from $-n$ to $n$. Hence for even $n$, they are $n+1$ consecutive even integers, 
while for odd $n$ they are $n+1$ consecutive odd integers.

\begin{rema}
	The eigenvectors of the matrix $C_n$ are also known, they can be expressed in terms of the Krawtchouk orthogonal polynomials \cite{Nom}.
\end{rema}

As the eigenvalues of \eqref{KSC} are known explicitly and because of the elegant and simple structure of both matrix and spectrum, $C_{n}$ is a standard test matrix for numerical eigenvalue computations (see e.g.~example~7.10 in \cite{GregoryKarney}), and part of some standard test matrix toolboxes (e.g.~\cite{Higham}). In MATLAB, $C_n$ can be produced from the gallery of test matrices using {\tt gallery(`clement',n+1)}. The Clement matrix also appears in several applications, e.g.~as the transition matrix of a Markov chain \cite{Abdel}, or in biogeography \cite{Igelnik}.

General tridiagonal matrices with closed form eigenvalues are rare, most examples being just variations of the tridiagonal matrix with fixed constants $a$, $b$ and $c$ on the subdiagonal, diagonal and superdiagonal respectively \cite{CuminatoMcKee,Yueh}. In this paper we present appealing two-parameter extensions of $C_{n}$ with closed form eigenvalues. These extensions first appeared in a different form in the paper \cite{OVdJ} as a special case of a class of matrices related to orthogonal polynomials. Special cases of these matrices were originally encountered in the context of finite quantum oscillator models (e.g.~\cite{JSV}) and their classification led to the construction of new interesting models \cite{OVdJ2}. 
Here, we feature them in a simpler, more accessible form which immediately illustrates their relation with $C_{n}$. Moreover, we consider some specific parameter values which yield interesting special cases.

Another purpose of this paper is to demonstrate by means of some numerical experiments the use of these extensions of $C_n$ as test matrices for numerical eigenvalue computations. 
Hereto, we examine how accurately the inherent MATLAB function {\tt eig()} is able to compute the eigenvalues of our test matrices compared to the exact known eigenvalues. An interesting feature of the new class of test matrices is that they include matrices with double eigenvalues for specific parameter values.

In section~\ref{sec:newmat} we present in an accessible form the two-parameter extensions of the Clement matrix. We state the explicit and rather simple form of their eigenvalues which makes them potentially interesting examples of eigenvalue test matrices. 
In section~\ref{sec:spec} we consider some specific parameter values for the new classes of test matrices which yield interesting special cases. 
In section~\ref{sec:num} we display some numerical results regarding the use of these extensions as test matrices for numerical eigenvalue computations. This is done by looking at the relative error when the exact known eigenvalues are compared with those computed using the inherent MATLAB function {\tt eig()}.

\section{New test matrices}
\label{sec:newmat}

Now, we consider the following extension of the matrix \eqref{KSC}, by generalizing its entries \eqref{entries} to:
\begin{equation}
\label{entriesex}
h_{k,k+1}=  \begin{cases} k  &\mbox{if $k$ even} \\ k+ a &\mbox{if $k$ odd}\end{cases}
\quad\mbox{ and }\quad
h_{n+2-k,n+1-k} = \begin{cases} k  &\mbox{if $k$ even} \\ k+ b &\mbox{if $k$ odd}\end{cases}
\end{equation}
where we introduce two parameters $a$ and $b$ 
(having a priori no restrictions). We will denote this extension by $H_{n}(a,b)$. 
For $n$ even, the matrix $H_{n}(a,b)$ is given by \eqref{Hahneven}. 
Note in particular that the first entry of the second row is $h_{2,1}=n$ and contains no parameter, and the same for $h_{n,n+1}=n$. For odd $n$, the matrix \eqref{Hahnodd} has entries $h_{2,1}=n+b$ and $h_{n,n+1}=n+a$ which now do contain parameters, contrary to the even case. 

The reason for considering this extension is that, similar to \eqref{KacEig} for $C_n$, we also have an explicit expression for the spectrum of $H_{n}(a,b)$, namely:

\begin{theo}\label{eve}
	For $n$ even, say $n=2m$, the $2m+1$ eigenvalues of 
	\begin{equation}
	H_{n}(a,b)= 
	\left(\begin{array}{ccccccc}
	0 & 1+a &  &   & &\\
	n & 0 &  2& & &\\
	& n-1+b& 0  & 3+a & &\\
	&   & \ddots &\ddots &\ddots  &  \\
	&   &  &   3+b  & 0 & n-1+a&  \\
	&   &      & & 2 & 0& n \\
	&   &  & & & 1+b& 0
	\end{array}\right)
	\label{Hahneven}
	\end{equation}
 are given by
	\begin{equation}
	\label{eigevenab}
	0, \pm\sqrt{(2k)(2k+a+b)}\qquad \mbox{for }k\in \{1,\dots,m\}.
	\end{equation}
\end{theo}

\begin{theo}\label{evo}
	For $n$ odd, say $n=2m+1$, the $2m+2$ eigenvalues of  
		\begin{equation}
		H_{n}(a,b)= 
		\left(\begin{array}{ccccccc}
		0 & 1+a &  &   & &\\
		n+b & 0 &  2& & &\\
		& n-1& 0  & 3+a & &\\
		&   & \ddots &\ddots &\ddots  &  \\
		&   &  &   3+b  & 0 & n-1&  \\
		&   &      & & 2 & 0& n+a \\
		&   &  & & & 1+b& 0
		\end{array}\right)
		\label{Hahnodd}
		\end{equation}
	are given by
	\begin{equation}
	\label{eigoddab}
	\pm \sqrt{(2k+1+a)(2k+1+b)}\qquad \mbox{for }k\in \{0,\dots,m\}
	\end{equation}
\end{theo}

We will prove the results for the symmetrized form of these matrices. We briefly elaborate on this. Consider a $(n+1)\times(n+1)$ tridiagonal matrix with zero diagonal
\begin{equation}
\label{A}
A= \left( \begin{array}{ccccc}
0 & b_1  &    0   &        &      \\
c_1 &  0  &  b_2  & \ddots &      \\
0  & c_2  &   0  & \ddots &  0   \\
&\ddots & \ddots & \ddots & b_{n} \\
&       &    0   &  c_{n}  &  0
\end{array} \right).
\end{equation} 
It is clear that the characteristic polynomial of $A$ depends on the products $b_ic_i$ ($i=1,\ldots,n$) only, and not on $b_i$ and $c_i$ separately. 
Therefore, if all the products $b_ic_i$ are positive, the eigenvalues of $A$ or of its symmetrized form
\begin{equation}
\label{A'}
A'= \left( \begin{array}{ccccc}
0 & \sqrt{b_1c_1}  &    0   &        &      \\
\sqrt{b_1c_1} &  0  &  \sqrt{b_2c_2}  & \ddots &      \\
0  & \sqrt{b_2c_2}   &   0  & \ddots &  0   \\
&\ddots & \ddots & \ddots & \sqrt{b_{n}c_n} \\
&       &    0   & \sqrt{b_{n}c_n}   &  0
\end{array} \right)
\end{equation} 
are the same. The eigenvectors of $A'$ are those of $A$ after multiplication by a diagonal matrix (the diagonal matrix that is used
in the similarity transformation from $A$ to $A'$).

Using this procedure, the aforementioned matrices can be made symmetric. For $C_n$ the entries \eqref{entries} can be symmetrized to
\[
\tilde{c}_{k,k+1} =\tilde{c}_{k+1,k} = \sqrt{k( n+1-k)} \qquad\mbox{for }k\in \{1,\dots,n\}.
\]
This matrix is also implemented in MATLAB, namely as {\tt gallery('clement',n+1,1)}.
For the extension $H_{n}(a,b)$, the entries of its symmetric form $\tilde H_{n}(a,b)$ are 
\[
\tilde{h}_{k,k+1}^e =\tilde{h}_{k+1,k}^e =   \begin{cases} \sqrt{k( n+1-k+b)}   &\mbox{if $k$ even} \\ \sqrt{(k+a)( n+1-k)}  &\mbox{if $k$ odd}\end{cases}
\]
for $n$ even and $k\in \{1,\dots,n\}$, while for $n$ odd we have
\[
\tilde{h}_{k,k+1}^o =\tilde{h}_{k+1,k}^o =  \begin{cases} \sqrt{k( n+1-k)}   &\mbox{if $k$ even} \\ \sqrt{(k+a)( n+1-k+b)} &\mbox{if $k$ odd}\end{cases}
\]
for $k\in \{1,\dots,n\}$. 

The above theorems are now proved using a property of the dual Hahn polynomials, which are defined in terms of the generalized hypergeometric series as follows \cite{Koekoek}
\begin{equation}\label{DHahn}
	R_n(\lambda(x);\gamma,\delta,N) =   
	{\;}_3F_2 \left( \genfrac..{0pt}{}{-x,x+\gamma+\delta+1,-n}{\gamma+1,-N}
	 ; 1 \right).
\end{equation}
The dual Hahn polynomials satisfy a discrete orthogonality relation, see \cite[(2.7)]{OVdJ}, and we denote the related orthonormal functions as $\tilde R_n(\lambda(x);\gamma,\delta,N)$. 
\begin{lemma}
The orthonormal dual Hahn functions satisfy the following pairs of recurrence relations:
\begin{align}
\notag	& \sqrt{(n+1+\gamma)(N-n)} \tilde R_n(\lambda(x);\gamma,\delta,N)
 - \sqrt{(n+1)(N-n+\delta)}\tilde R_{n+1}(\lambda(x);\gamma,\delta,N)  \\ \label{lem1} &=  \sqrt{x(x+\gamma+\delta+1)} \tilde R_n(\lambda(x-1);\gamma+1,\delta+1,N-1), \\
\notag	& -\sqrt{(n+1)(N-n+\delta)}\tilde R_n(\lambda(x-1);\gamma+1,\delta+1,N-1)
 +\sqrt{(n+2+\gamma)(N-n-1)}\\ \label{lem2} & \times \tilde R_{n+1}(\lambda(x-1);\gamma+1,\delta+1,N-1) = \sqrt{x(x+\gamma+\delta+1)} \tilde R_{n+1}(\lambda(x);\gamma,\delta,N)
\end{align}
and
\begin{align}
\notag &\sqrt{(n+1+\gamma)(N-n+\delta)}\tilde R_n(\lambda(x);\gamma,\delta,N)
-\sqrt{(n+1)(N-n)}  \tilde R_{n+1}(\lambda(x);\gamma,\delta,N) \\ \label{lem3} &=\sqrt{(x+\gamma+1)(x+\delta)}  \tilde R_n(\lambda(x);\gamma+1,\delta-1,N),\\
\notag	&
-\sqrt{(n+1)(N-n)} \tilde R_n(\lambda(x);\gamma+1,\delta-1,N)
+ \sqrt{(n+2+\gamma)(N-n+\delta-1)} \\ \label{lem4} & \times  \tilde R_{n+1}(\lambda(x);\gamma+1,\delta-1,N)  = \sqrt{(x+\gamma+1)(x+\delta)} \tilde R_{n+1}(\lambda(x);\gamma,\delta,N).
\end{align}
\end{lemma}
\noindent  {\bf Proof.}

The first two relations follow from the case dual Hahn I of \cite[Theorem 1]{OVdJ} multiplied by the square root of the weight function and norm squared, 
and similarly the last two from the case dual Hahn III.

\hfill$\Box$

\noindent  {\bf Proof of Theorem 2.}

Let $n$ be an even integer, say $n=2m$, and let $a$ and $b$ be real numbers greater than $-1$.  
Take $k\in \{1,\dots,m\}$ and let $U_{\pm k}=(u_1,\dotsc,u_{n+1})^T$ be the column vector with entries
\[
u_l = \begin{cases} 
(-1)^{(l-1)/2} \tilde R_{(l-1)/2}(\lambda(k);\frac{a-1}{2},\frac{b-1}{2},m)  &\mbox{if $l$ odd} 
\\
 \pm (-1)^{l/2-1} \tilde R_{l/2-1}(\lambda(k-1);\frac{a+1}{2},\frac{b+1}{2},m-1) &\mbox{if $l$ even.}
\end{cases}
\]
We calculate the entries of the vector $\tilde H_{n}(a,b)U_{\pm k}$ to be 
\[
(\tilde H_{n}(a,b)U_{\pm k})_l = \tilde{h}_{l,l-1}^e u_{l-1} + \tilde{h}_{l,l+1}^e u_{l+1}.
\]
For $l$ even, using the recurrence relation \eqref{lem1} with the appropriate parameter values substituted in the orthonormal dual Hahn functions, this becomes
\begin{align*}
(\tilde H_{n}(a,b)U_{\pm k})_l = \ & \sqrt{(l-1+a)( 2m+2-l)}  (-1)^{l/2-1} \tilde R_{l/2-1}(\lambda(k);\tfrac{a-1}{2},\tfrac{b-1}{2},m)    
\\ & +\sqrt{l( 2m+1-l+b)} (-1)^{l/2} \tilde R_{l/2}(\lambda(k);\tfrac{a-1}{2},\tfrac{b-1}{2},m) 
\\  = \ &   2\sqrt{k(k+\tfrac{a}{2}+\tfrac{b}{2})} (-1)^{l/2-1} \tilde R_{l/2-1}(\lambda(k-1);\tfrac{a+1}{2},\tfrac{b+1}{2},m-1) \\  = \ & \pm \sqrt{(2k)(2k+a+b)} \ u_l.
\end{align*}
Similarly, for $l$ odd we have, using now the recurrence relation \eqref{lem2}, 
\begin{align*}
(\tilde H_{n}(a,b)U_{\pm k})_l = \ & \pm\sqrt{(l-1)( 2m+2-l+b)} (-1)^{(l-3)/2} \tilde R_{(l-3)/2}(\lambda(k-1);\tfrac{a+1}{2},\tfrac{b+1}{2},m-1)    \\ & 
\pm\sqrt{(l+a)( 2m+1-l)}  (-1)^{(l-1)/2} \tilde R_{(l-1)/2}(\lambda(k-1);\tfrac{a+1}{2},\tfrac{b+1}{2},m-1) \\  = \ &   \pm2\sqrt{k(k+\tfrac{a}{2}+\tfrac{b}{2})} (-1)^{(l-1)/2} \tilde R_{(l-1)/2}(\lambda(k);\tfrac{a-1}{2},\tfrac{b-1}{2},m) \\  = \ & \pm \sqrt{(2k)(2k+a+b)} \ u_l.
\end{align*}
Finally, define $U_{0}=(u_1,\dotsc,u_{n+1})^T$ as the column vector with entries
\[
u_l = \begin{cases} 
(-1)^{(l-1)/2} \tilde R_{(l-1)/2}(\lambda(0);\frac{a-1}{2},\frac{b-1}{2},m)  &\mbox{if $l$ odd} 
\\
0 &\mbox{if $l$ even.}
\end{cases}
\]
Putting $x=0$ on the right-hand side of \eqref{lem1}, it is clear that the entries of the vector $\tilde H_{n}(a,b)U_{0}$ are all zero.

This shows that the eigenvalues of $\tilde H_{n}(a,b)$ are given by \eqref{eigevenab}, so its characteristic polynomial must be
\[
\lambda \prod_{k=1}^m \bigl(\lambda^2- (2k)(2k+a+b)\bigr),
\]
which allows us to extend the result to arbitrary parameters $a$ and $b$.

\hfill$\Box$

Theorem 3 is proved in the same way, using now relations \eqref{lem3} and \eqref{lem4}.

\section{Special cases}
\label{sec:spec}

We now consider some particular cases where the eigenvalues as given in
Theorem~\ref{eve} and Theorem~\ref{evo} reduce to integers or have a special form.

For the specific values $a=0$ and $b=0$, it is clear that $H_{n}(0,0)$ reduces to $C_n$ for both even and odd values of $n$. As expected, the explicit formulas for the eigenvalues \eqref{eigevenab} and \eqref{eigoddab} also reduce to \eqref{KacEig}.

Next, we look at $n$ even. In \eqref{eigevenab}, we see that the square roots cancel if we take $b=-a$. For this choice of parameters, the eigenvalues of $H_{n}(a,-a)$ are even integers, which are precisely the same eigenvalues as those of the Clement matrix \eqref{KacEig}. However, the matrix 
\begin{equation}
H_{n}(a,-a)= 
\left(\begin{array}{ccccccc}
0 & 1+a &  &   & &\\
n & 0 &  2& & &\\
& n-1-a& 0  & 3+a & &\\
&   & \ddots &\ddots &\ddots  &  \\
&   &  &   3-a  & 0 & n-1+a&  \\
&   &      & & 2 & 0& n \\
&   &  & & & 1-a& 0
\end{array}\right)
\label{Hahnevena}
\end{equation}
still contains a parameter $a$ which does not affect its eigenvalues. 
So for every even integer $n$, \eqref{Hahnevena} gives rise to a one-parameter family of tridiagonal matrices with zero diagonal whose eigenvalues are given by \eqref{KacEig}. This property is used explicitly in \cite{OVdJ2} to construct a finite oscillator model with equidistant position spectrum.
 
For $n$ odd, say $n=2m+1$, the square roots in \eqref{eigoddab} cancel if we take $b=a$. 
Substituting $a$ for $b$, the matrix \eqref{Hahnodd} becomes
\begin{equation}
H_{n}(a,a)= 
\left(\begin{array}{ccccccc}
0 & 1+a &  &   & &\\
n+a & 0 &  2& & &\\
& n-1& 0  & 3+a & &\\
&   & \ddots &\ddots &\ddots  &  \\
&   &  &   3+a  & 0 & n-1&  \\
&   &      & & 2 & 0& n+a \\
&   &  & & & 1+a& 0
\end{array}\right),
\label{Hahnodda}
\end{equation}
while for the eigenvalues we get 
\begin{equation}
\label{eigodda}
\pm \lvert 2k+1+a \rvert \qquad \mbox{for }k\in \{0,\dots,m\}.
\end{equation}
These are integers for integer $a$ and real numbers for real $a$. We see that for $a=0$, \eqref{eigodda} reduces to the eigenvalues \eqref{KacEig}, but then the matrix is precisely $C_n$. Non-zero values for $a$ induce a shift in the eigenvalues, away from zero for positive $a$ and towards zero for $-1<a<0$. However, for $-n<a<-1$ the positive and negative (when $a>-1$) eigenvalues get mingled. Moreover, for $a$ equal to a negative integer ranging from $-1$ to $-n$, we see that there are double eigenvalues. A maximum number of double eigenvalues occurs for $a=-m-1$, then each of the values
\[
2k-m\qquad \mbox{for }k\in \{0,\dots,m\},
\]
is a double eigenvalue. By choosing $a$ nearly equal to a negative integer, we can produce a matrix with nearly, but not exactly, equal eigenvalues. 
For $a<-n$, all positive eigenvalues (when $a>-1$) become negative and vice versa. Finally, for the special value $a=-n-1$, the eigenvalues \eqref{eigodda} also reduce to \eqref{KacEig} while the matrix becomes
\begin{equation*}
H_{n}(-n-1,-n-1)= 
\left(\begin{array}{ccccccccc}
0 & -n &  &   & &\\
-1 & 0 &  2& & &\\
& n-1& 0  & 2-n & &\\
& & -3& 0  & 4 & &\\
& &   & \ddots &\ddots &\ddots  &  \\
& &   &  &   4  & 0 & -3&  \\
&& &   &  &   2-n  & 0 & n-1&  \\
&& &   &      & & 2 & 0& -1 \\
&& &   &  & & & -n& 0
\end{array}\right).
\end{equation*}
This is up to a similarity transformation, as explained at the end of section~2, the matrix $C_n$.

Also for $n$ odd, another peculiar case occurs when $b=a=1$. Scaling by one half we then have the matrix
\begin{equation}
\frac{1}{2}H_{2m+1}(1,1)= 
\left(\begin{array}{cccccccccc}
0 & 1 &  &   & &\\
m+1 & 0 &  1& & &\\
& m& 0  & 2 & &\\
&& m& 0  & 2 & &\\
&   & & \ddots &\ddots &\ddots  &  \\
&   & & &   2  & 0 & m&  \\
&  & & & &   2  & 0 & m&  \\
& &  &  &    & & 1 & 0& m+1 \\
&  & &  && & & 1& 0
\end{array}\right).
\label{Hahnodd11}
\end{equation}
with eigenvalues, by \eqref{eigoddab},
\[
	\pm (k+1)\qquad \mbox{for }k\in \{0,\dots,m\}.
\]
The even equivalent of this matrix, 
\begin{equation} 
\frac{1}{2}H_{2m}(1,1)= \left(\begin{array}{cccccccccc}
0 & 1 &  &   & &\\
m & 0 &  1& & &\\
& m& 0  & 2 & &\\
&& m-1& 0  & 2 & &\\
&   & & \ddots &\ddots &\ddots  &  \\
&   & & &   2  & 0 & m-1&  \\
&  & & & &   2  & 0 & m&  \\
& &  &  &    & & 1 & 0& m \\
&  & &  && & & 1& 0
\end{array}\right),
\label{Hahneven11}
\end{equation}
does not have integer spectrum, but instead, using the expression \eqref{eigevenab}, has as eigenvalues
\[
0,\pm \sqrt{k(k+1)}\qquad \mbox{for }k\in \{1,\dots,m\}.
\]

We have reviewed the special cases where the explicit formulas for the eigenvalues \eqref{eigevenab} and \eqref{eigoddab} which generally contain square roots, reduce to integers. 
In the following, we will use the notation 
\begin{equation}
\label{speccase}
H_{n}(a)=\begin{cases} H_{n}(a,-a)  &\mbox{if $n$ even} \\ H_{n}(a,a)  &\mbox{if $n$ odd}\end{cases}
\end{equation}
to denote these special cases.

\section{Numerical results}
\label{sec:num}

We now examine some numerical results regarding the use of the extensions of the previous sections as test matrices for numerical eigenvalue computations. This is done by comparing the exact known eigenvalues of $H_n(a,b)$ with those computed using the inherent MATLAB function {\tt eig()}. These numerical experiments are included merely to illustrate the potential use of the matrices $H_n(a,b)$ as eigenvalue test matrices, and to examine the sensitivity of the computed eigenvalues on the new parameters.

A measure for the accuracy of the computed eigenvalues is the relative error
\[
\frac{\lVert x-x^* \rVert_{\infty}}{\lVert x\rVert_{\infty}},
\]
where $x$ is the vector of eigenvalues as ordered list (by real part) and $x^*$ its approximation.

Recall that both for $n$ odd and $n$ even, for the special case $H_{n}(a)$ the square roots in the expressions for the eigenvalues cancel, yielding real eigenvalues for every real value of the parameter $a$. 
In the general case, the eigenvalues \eqref{eigevenab} are real when $a+b>-2$ and those in \eqref{eigoddab} are real when $a>-1$ and $b>-1$ or $a<-n$ and $b<-n$. A first remark is that when we compute the spectrum of $C_n$ using {\tt eig()} in MATLAB, eigenvalues with imaginary parts are found when $n$ exceeds 116, but not for lower values of $n$. Therefore, for the extensions, we have chosen $n=100$ and $n=101$ (for the even and the odd case respectively) for most of our tests, as this gives reasonably large matrices but is below the bound of 116. We will see that in this case for the extensions, the {\tt eig()} function in MATLAB does find eigenvalues with imaginary part for certain parameter values.



\begin{figure}[!htbp]
	\begin{tabular}{cc}
		\hline \\
		\includegraphics[scale=0.525]{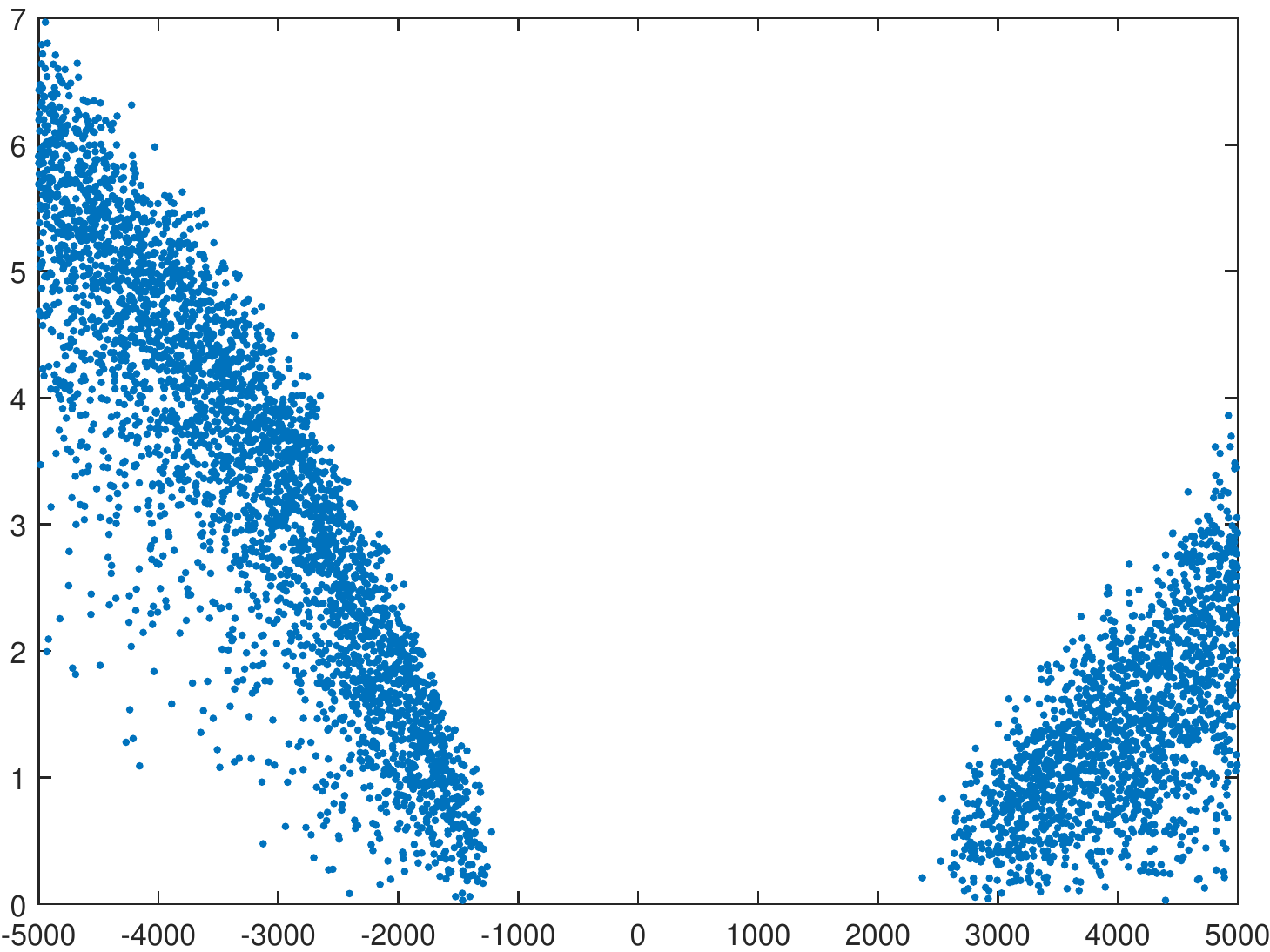} &
		\includegraphics[scale=0.525]{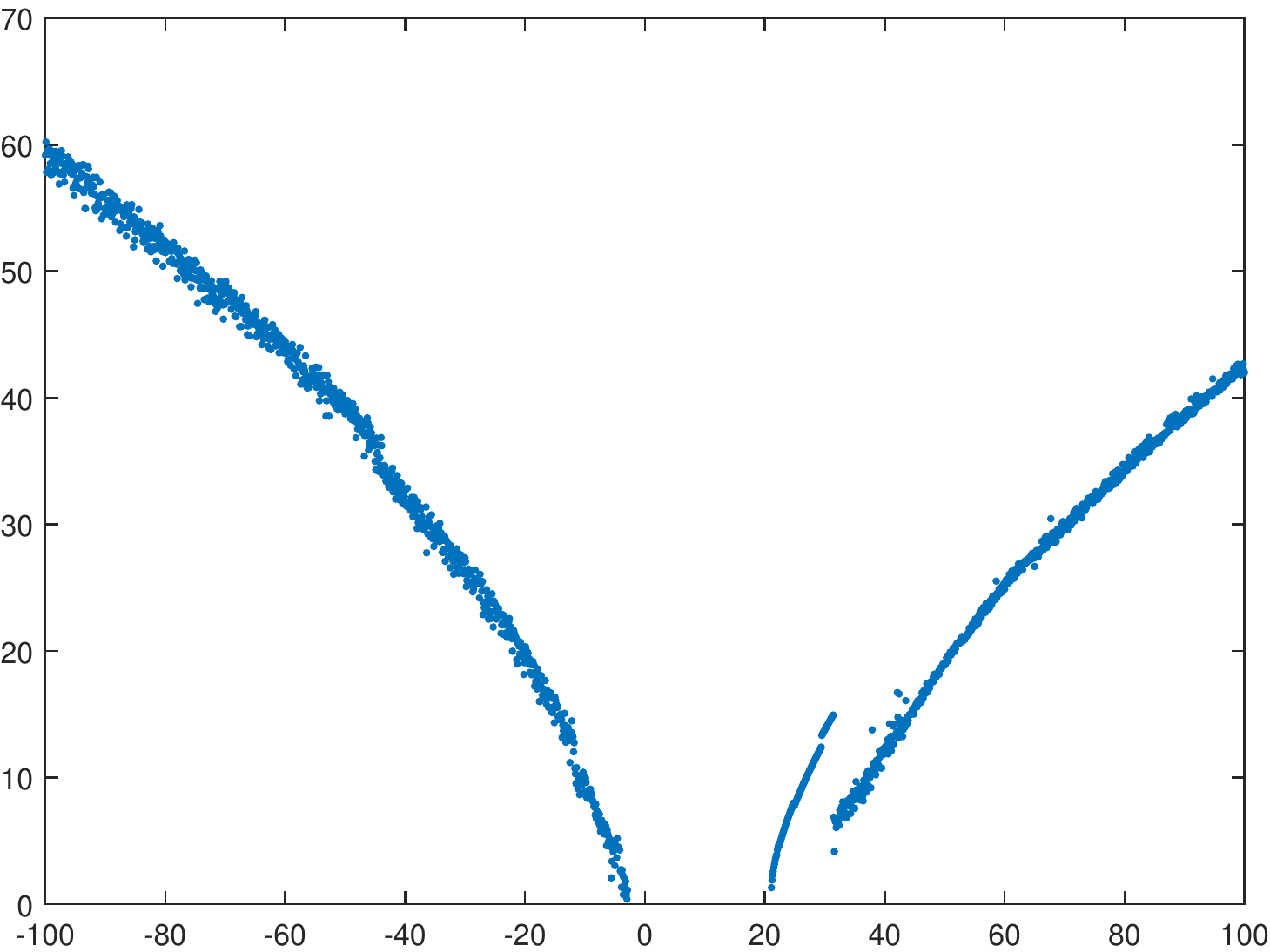}\\
		\\
		\hline
	\end{tabular} 
	\caption{Plots of the largest imaginary part of the computed eigenvalues of $H_n(a)$ ($n$ even and $a=-b$) for different values of $a$ on horizontal axis. On the left for $n=10$, on the right for $n=100$.}
	\label{fig1}
\end{figure}

\begin{figure}[!htbp]
	\begin{tabular}{cc}
		\hline \\
		\includegraphics[scale=0.525]{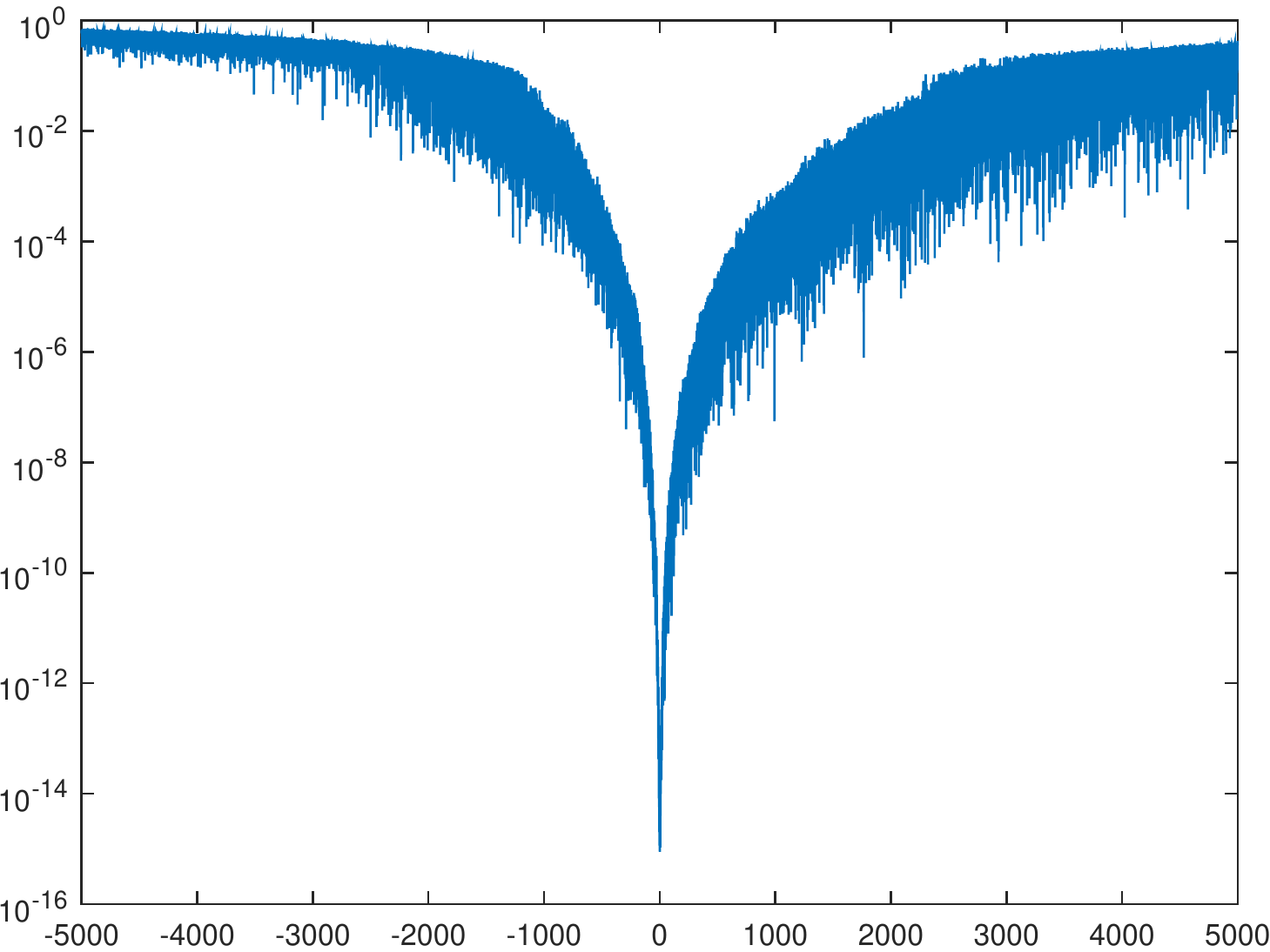} & \includegraphics[scale=0.525]{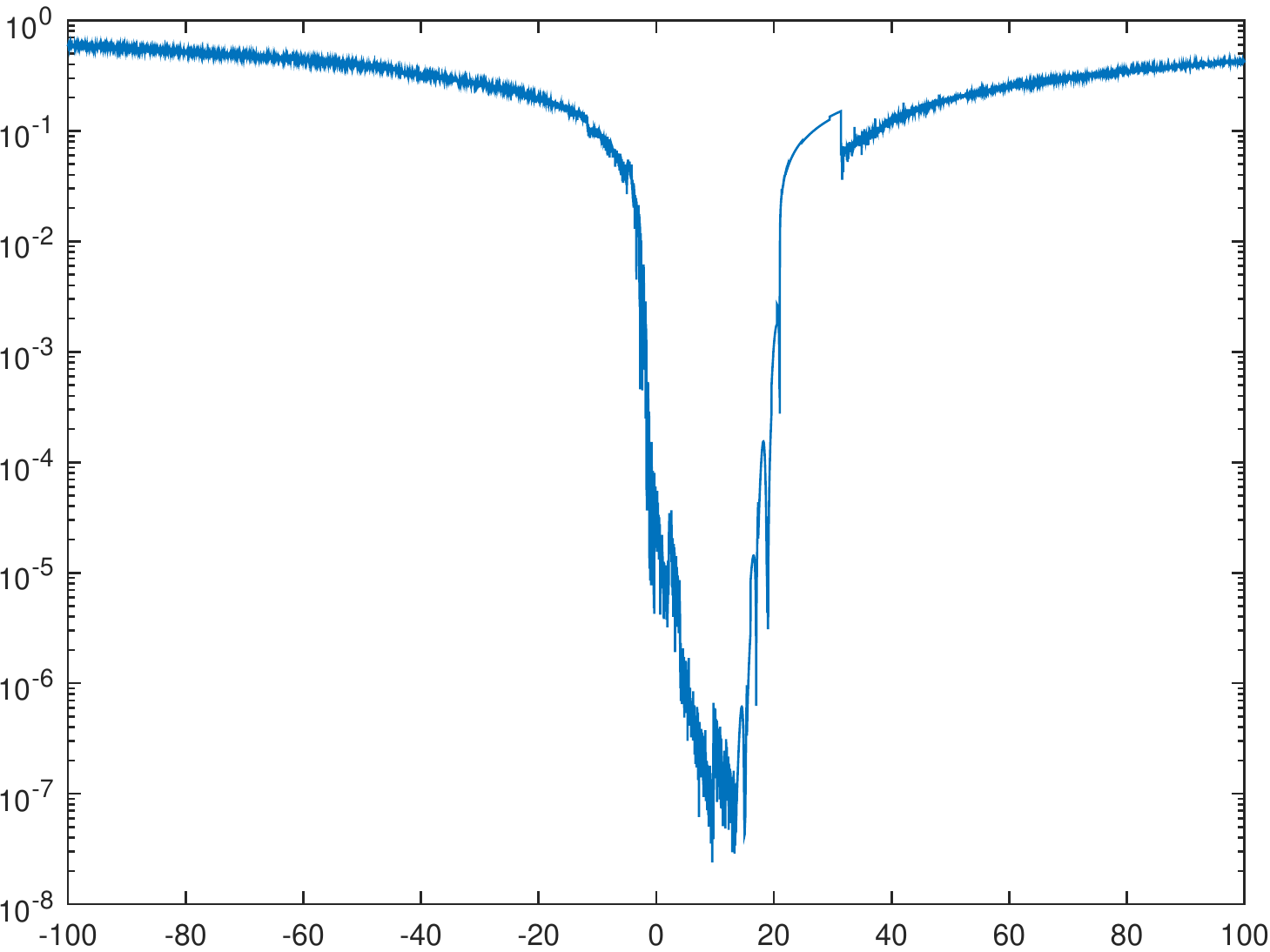} \\
		\hline
	\end{tabular} 
	\caption{Plots of the relative error in the computed eigenvalues of 
		$H_n(a)$
		($n$ even and $a=-b$) for different values of $a$ on horizontal axis. Left $n=10$ and right $n=100$.}
	\label{fig2}
\end{figure}

We first consider the special case \eqref{speccase}. For $n$ even, $H_{n}(a)$ has the eigenvalues \eqref{KacEig}, which are integers independent of the parameter $a$. 
In figure~\ref{fig1}, we have depicted the largest imaginary part of the computed eigenvalues for the matrix $H_{n}(a)$ for $n=10$ and $n=100$ at different values for the parameter $a$. We see that outside a central region imaginary parts are found. 
For example, for $H_{100}(a)$, MATLAB finds eigenvalues with imaginary parts when $a>21$ or $a<-2.5$.  
Moreover, the relative error for the computed eigenvalues, shown in figures~\ref{fig2}, increases as $a$ approaches the region where eigenvalues with imaginary parts are found. In this latter region, the size of the relative error is of course due to the presence of imaginary parts which do not occur in the theoretical exact expression for the eigenvalues.
As a reference, the relative error for $C_{100}$ is $3.6612 \times 10^{-5}$, while that for $H_{100}(20)$ is $1.1471 \times 10^{-3}$ and $4.9444 \times 10^{-3}$ for $H_{100}(20.97)$.

\begin{figure}[!htbp]
	\begin{tabular}{cc}
		\hline \\
		\includegraphics[scale=0.525]{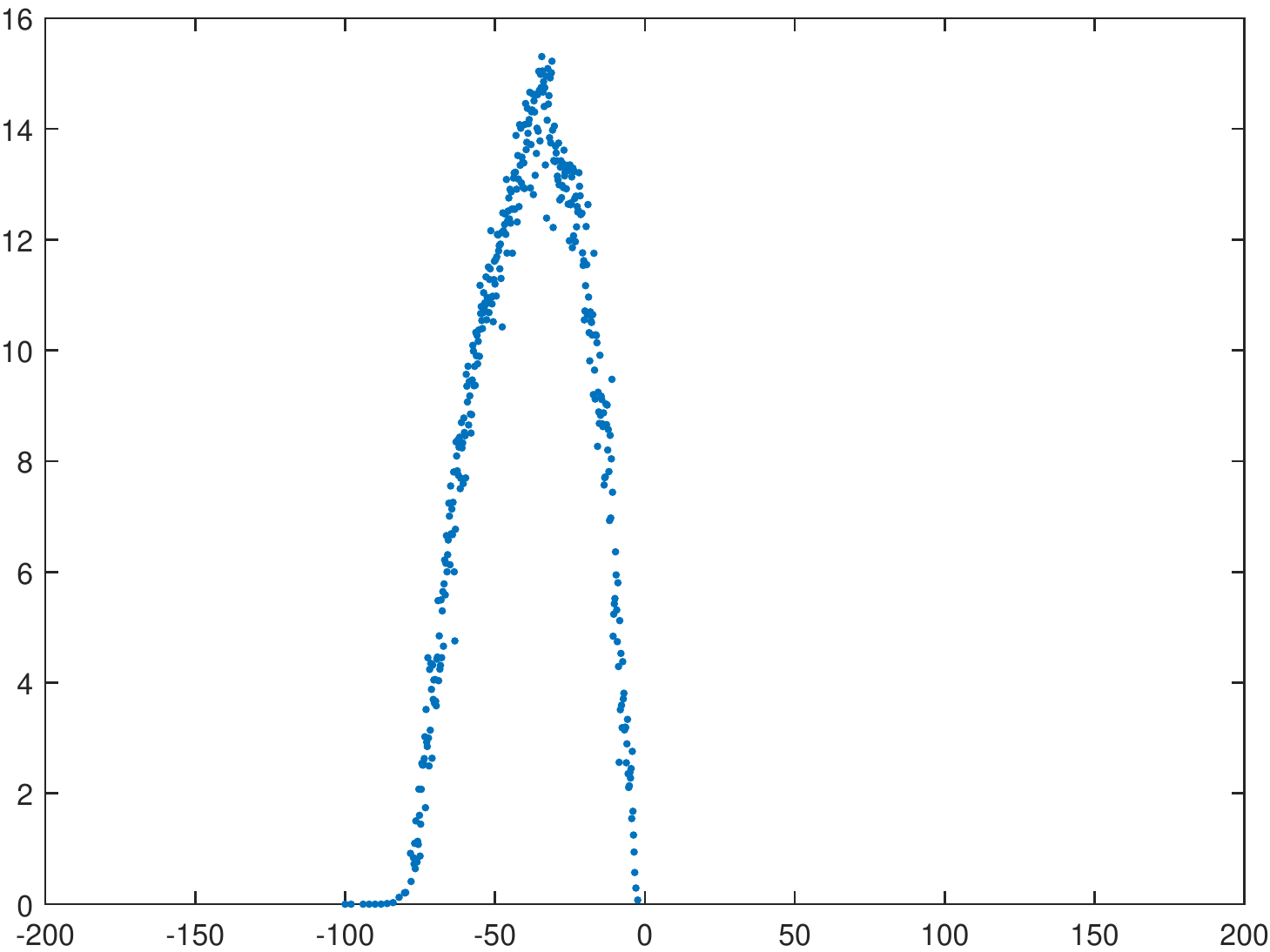} & 
		\includegraphics[scale=0.525]{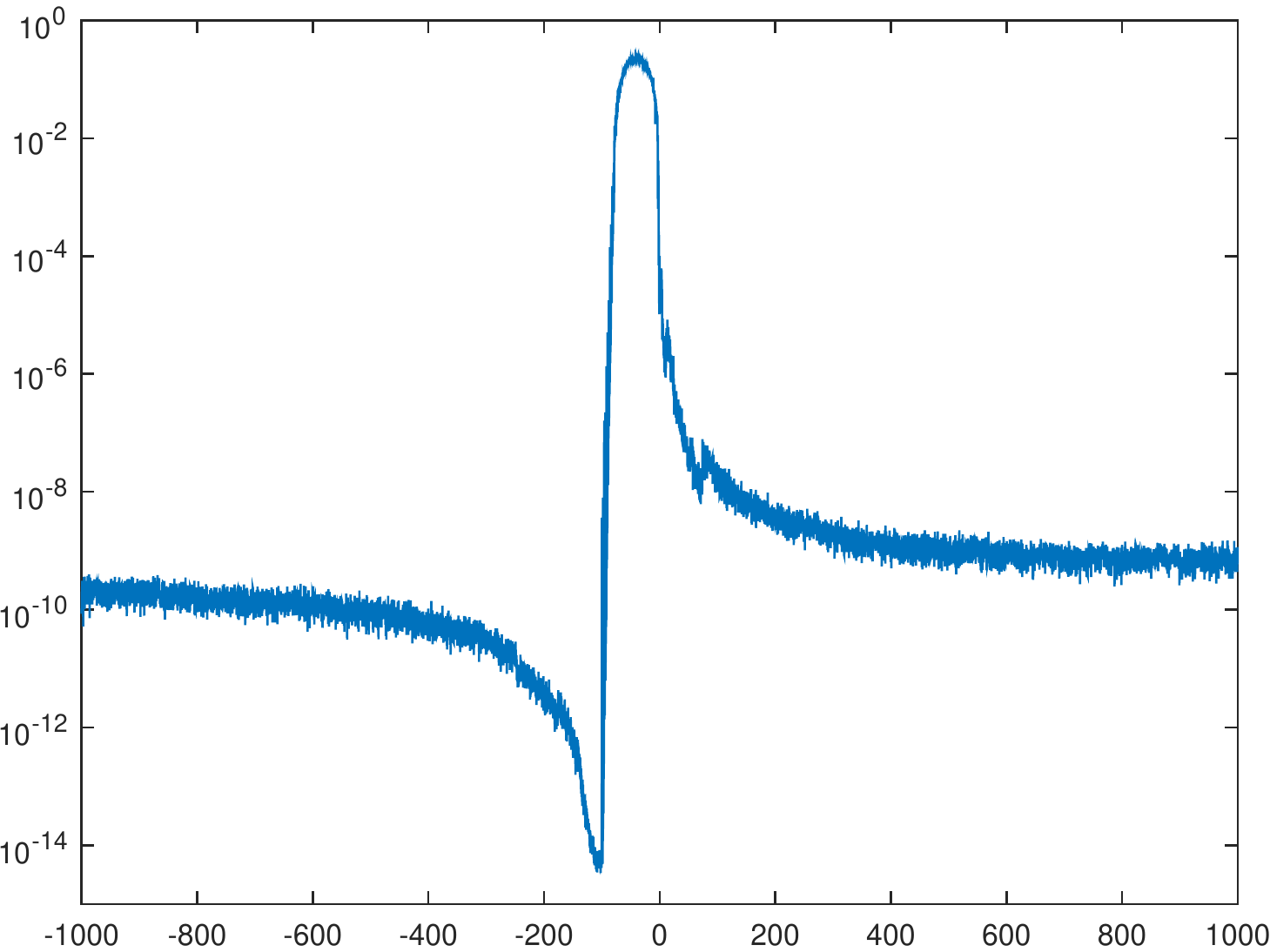}\\
		\\
		\hline
	\end{tabular} 
	\caption{For different values of $a$ as denoted on the horizontal axis, on the left a plot of the largest imaginary part and on the right a plot of the relative error of the computed eigenvalues of $H_{101}(a)$.}
	\label{fig3}
\end{figure}

For $n$ odd, $H_{n}(a)$ has the eigenvalues \eqref{eigodda}, which dependent on the parameter $a$ but are real for every real number $a$. Nevertheless, even for a small dimension such as $n=11$, eigenvalues with (small) imaginary parts are found when $a$ equals $-2,-4,-6$ or $-8$. This produces a relative error of order $10^{-8}$, while for other values of $a$ (and for the Clement matrix) the relative error is of order $10^{-15}$, near machine precision. 
For $H_{101}(a)$, the largest imaginary part of the computed eigenvalues is portrayed in figure~\ref{fig3}, together with the relative error. 
MATLAB finds eigenvalues with imaginary parts when $-100\leq a<-1.5$. These findings correspond to the region where double eigenvalues occur as mentioned in the previous section. The relative error is largest around this region and is several orders smaller when moving away from this region. As a reference, the relative error for $C_{101}$ is $3.6881 \times 10^{-5}$, while that for $H_{101}(-1.75)$ is $1.4840 \times 10^{-3}$. 
Finally, we note that eigenvalues with imaginary parts also appear when $a$ is extremely large, i.e.~$a>10^{10}$ or $a<-10^{10}$.



\begin{figure}[!htbp]
	\begin{tabular}{ccc}
		\hline \\
		\includegraphics[scale=0.525]{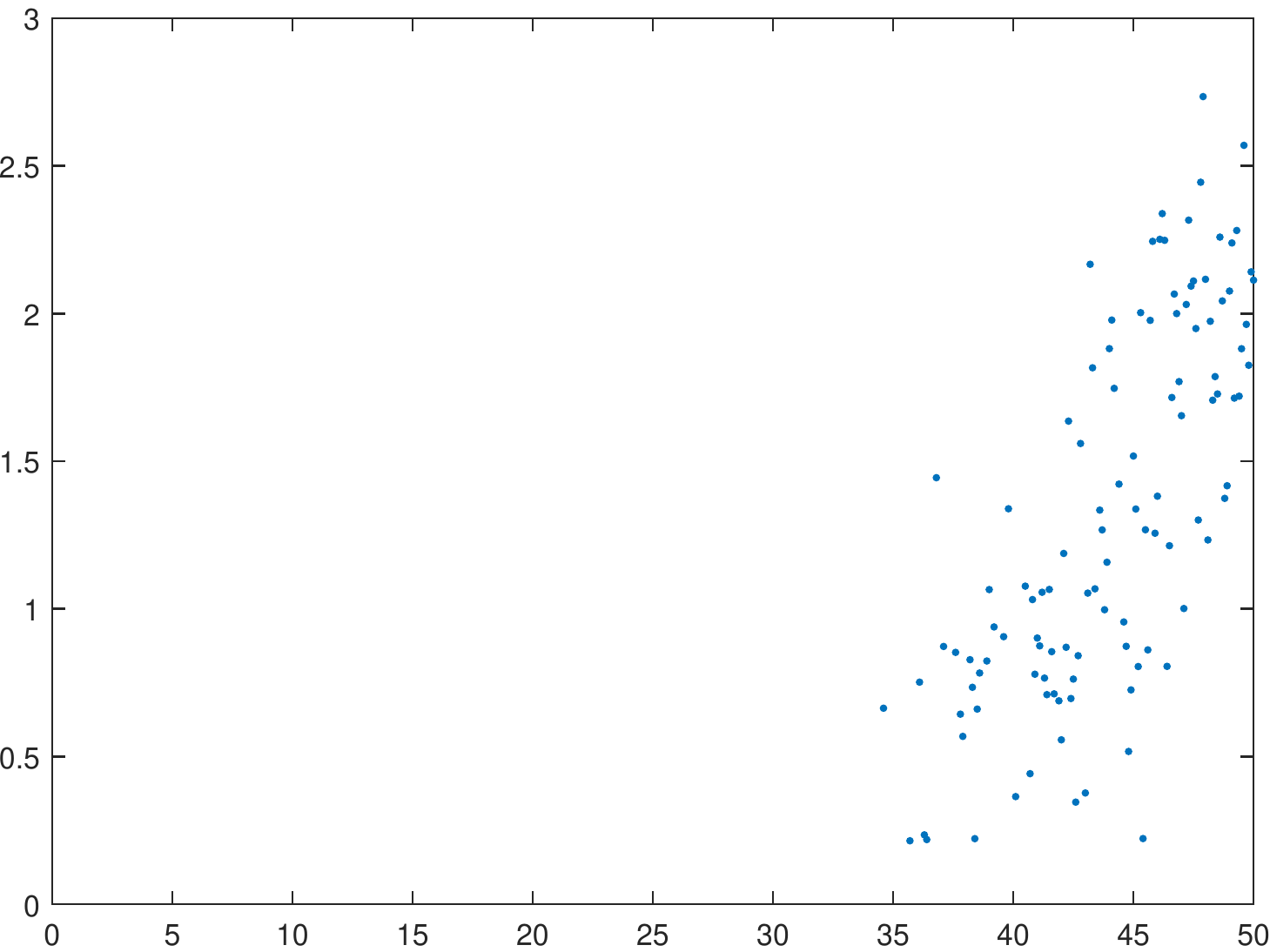}  &  \includegraphics[scale=0.525]{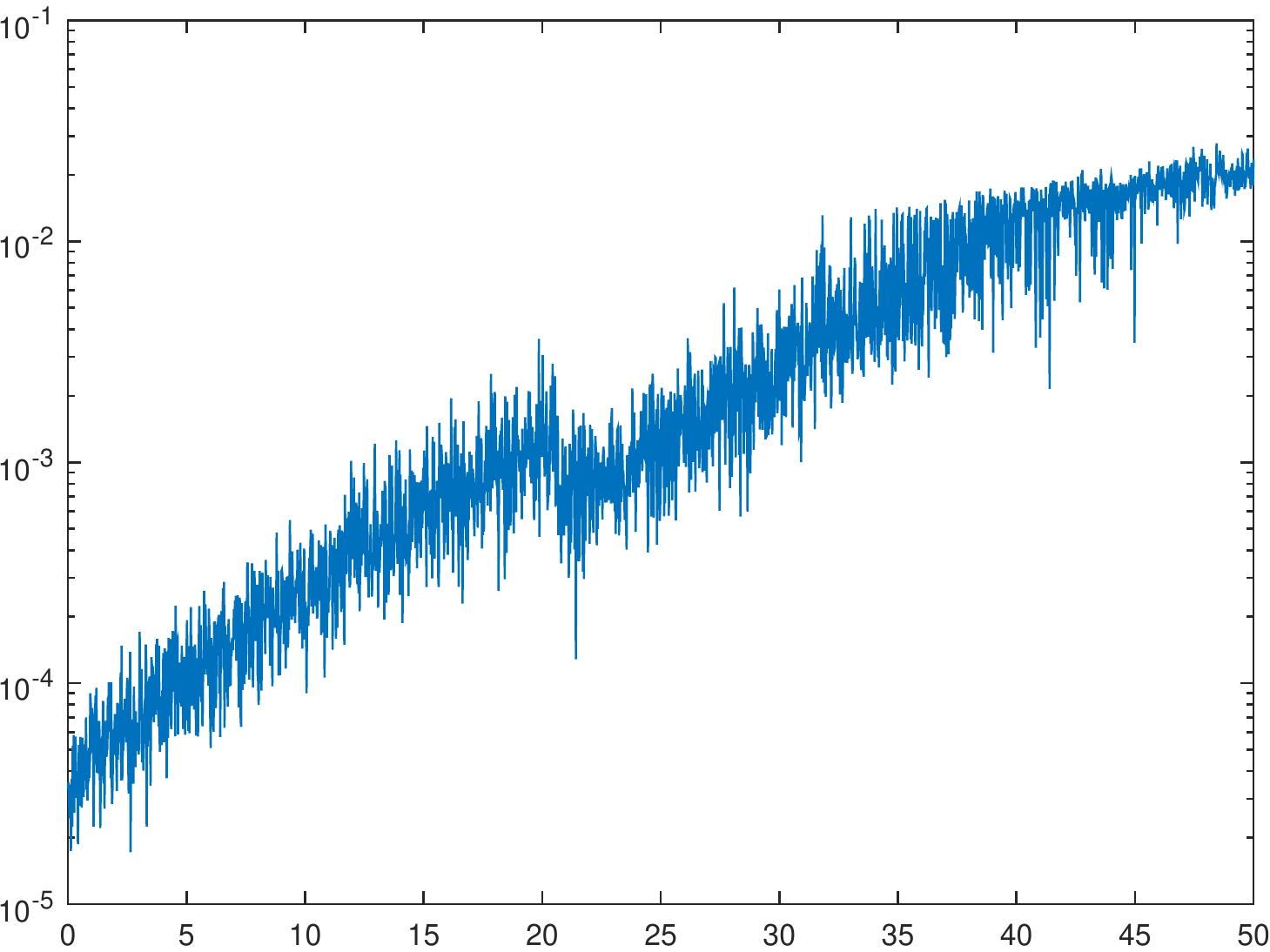} \\
		\includegraphics[scale=0.525]{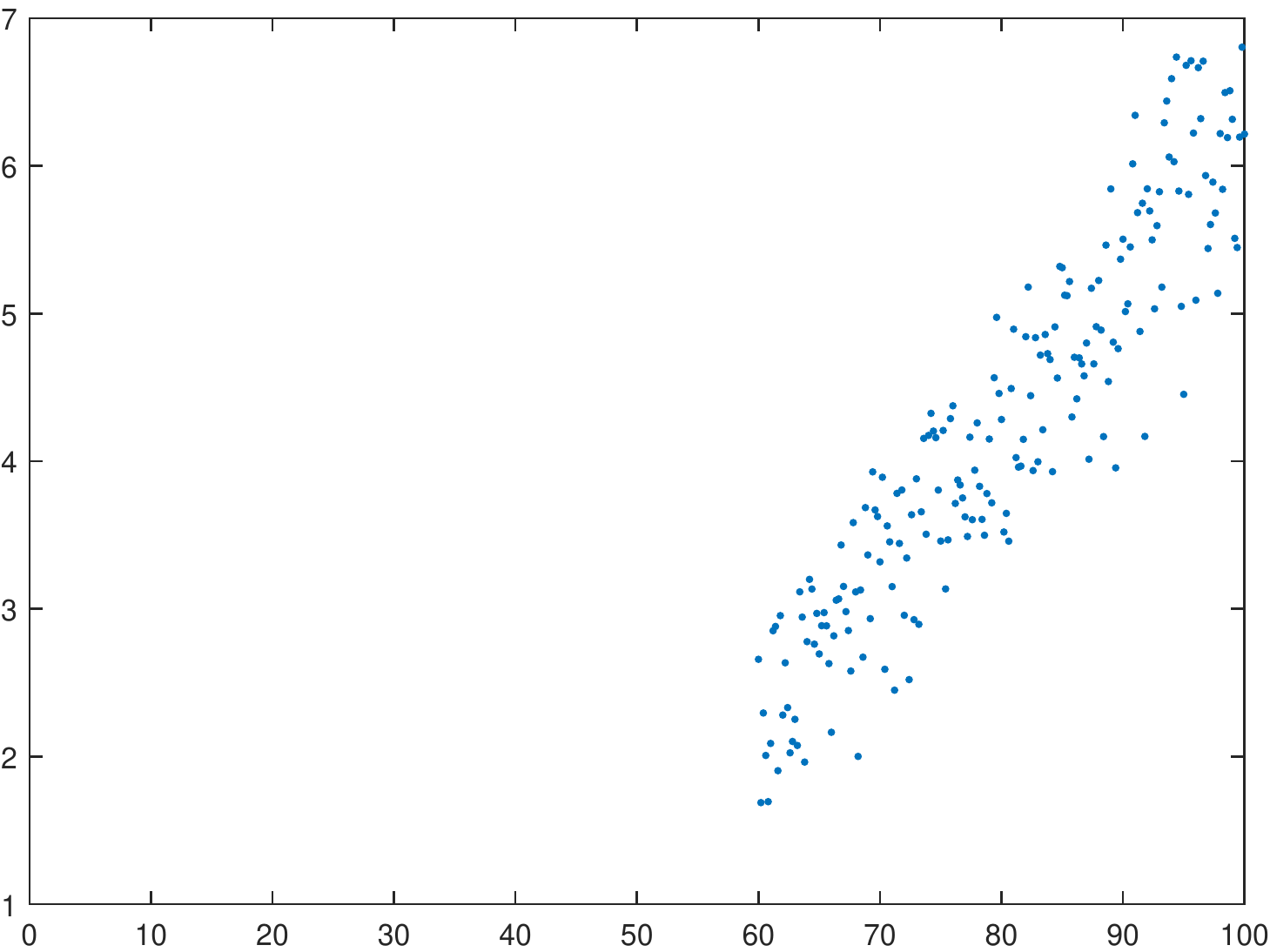} & \includegraphics[scale=0.525]{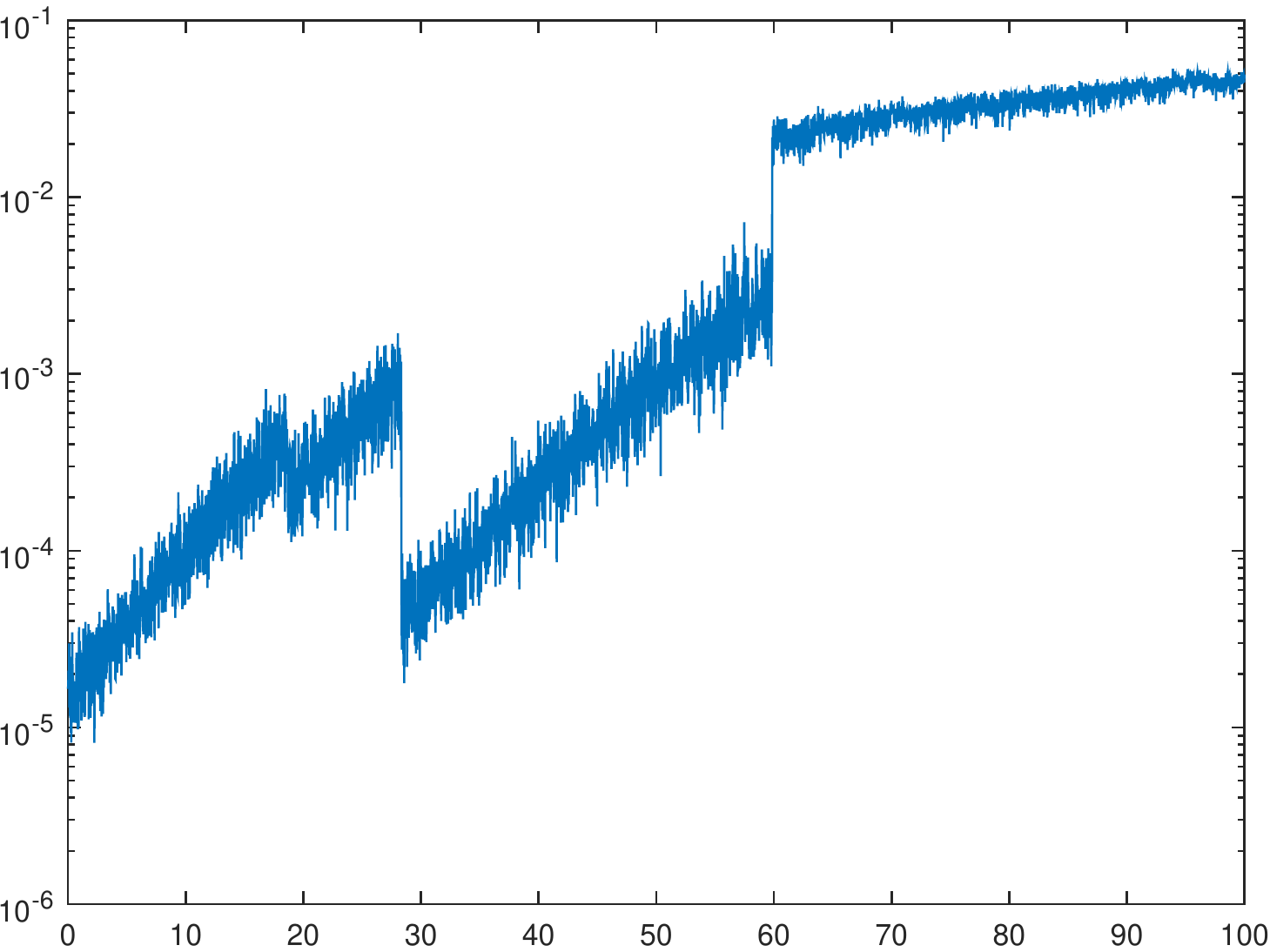}  \\
		\includegraphics[scale=0.525]{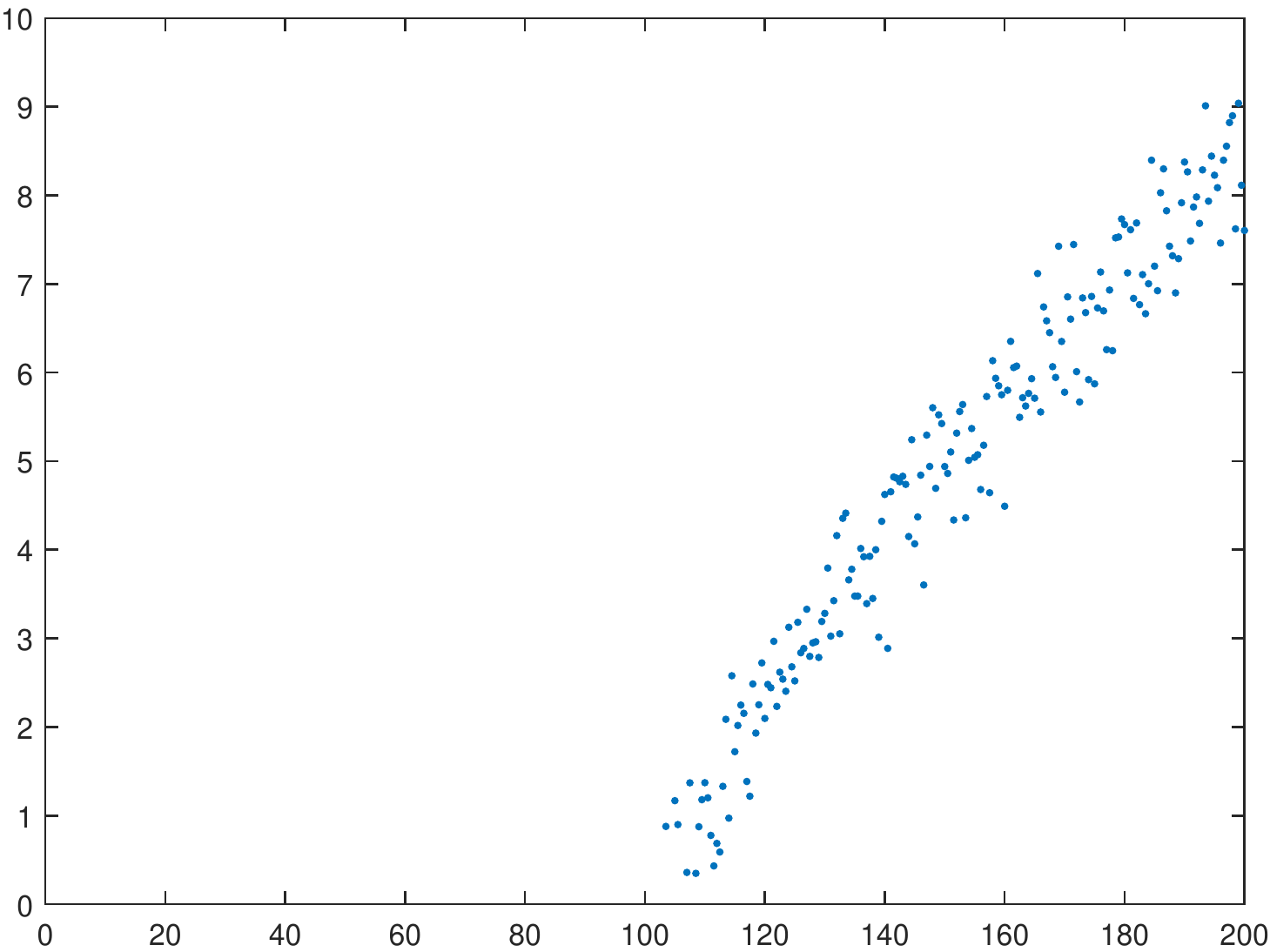} & \includegraphics[scale=0.525]{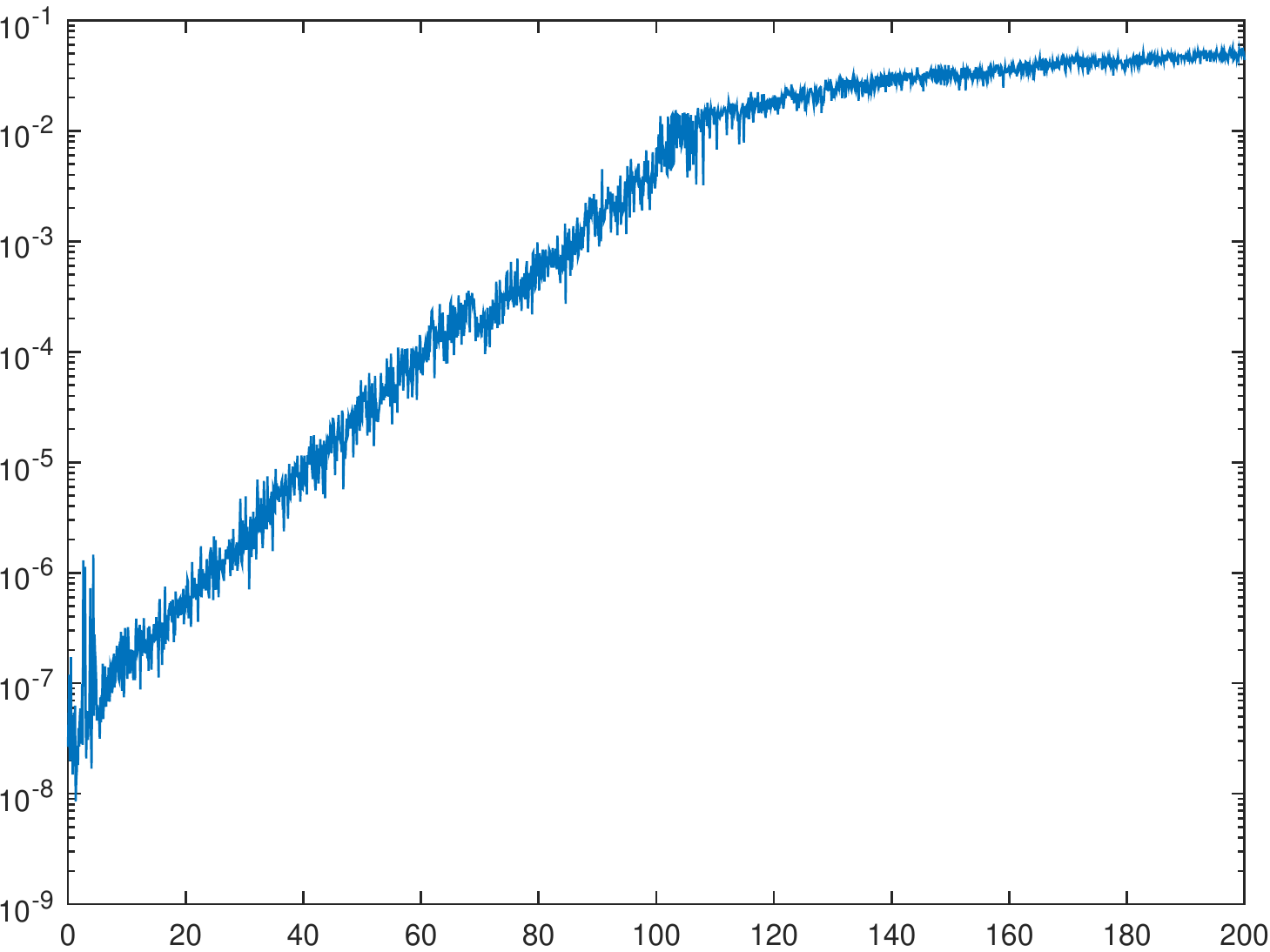}\\
		\\
		\hline
	\end{tabular} 
	\caption{Plots of the largest imaginary part (left) and the relative error (right) in the computed eigenvalues of $H_{100}(a,b)$, horizontal axis varying values of $b$. Top row $a=0$, middle row $a=1$, bottom row $a=20$.}
	\label{fig4}
\end{figure}

\begin{figure}[!htbp]	
	\begin{tabular}{ccc}
		\hline \\
		\includegraphics[scale=0.525]{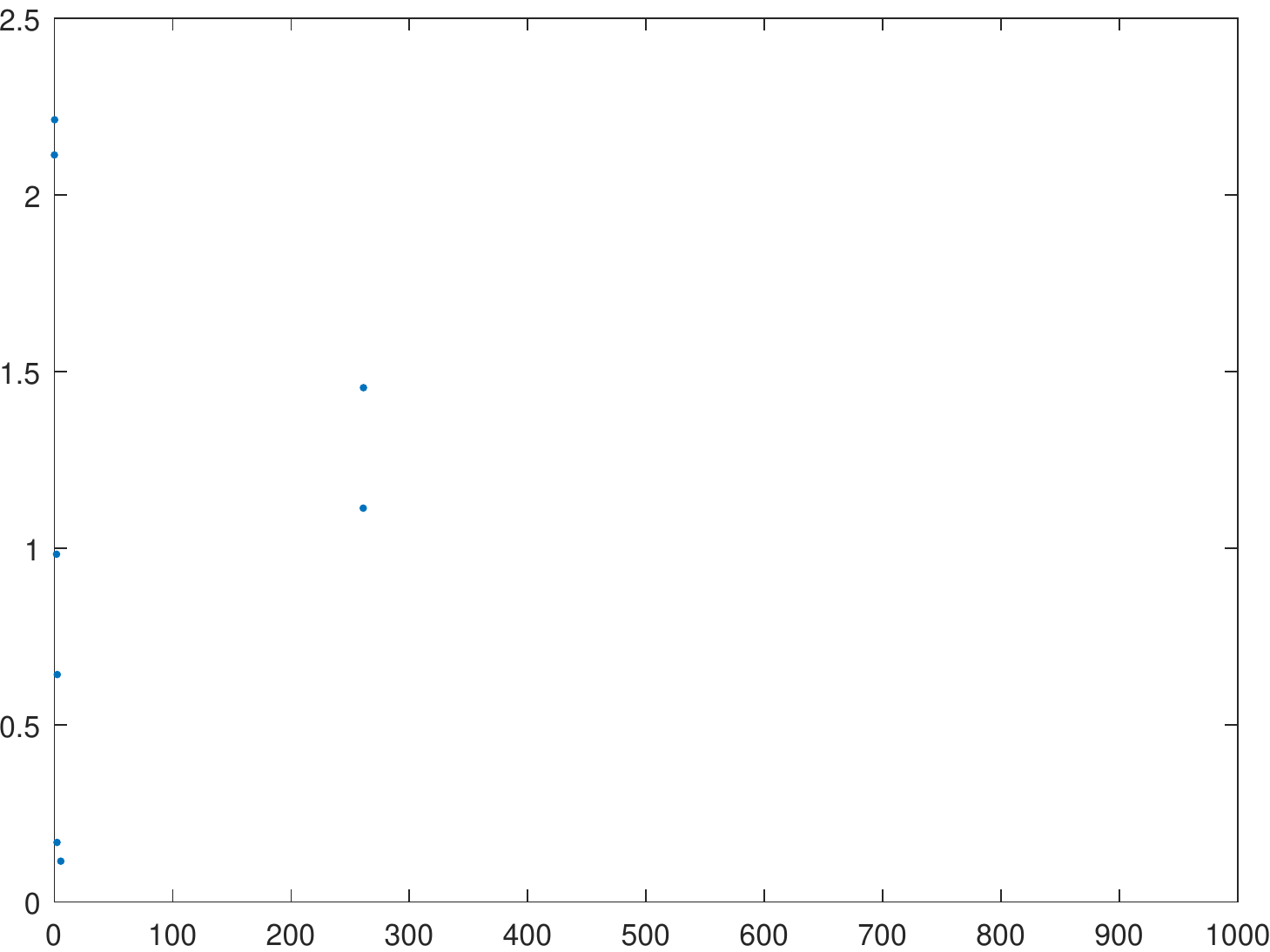}  & 
		\includegraphics[scale=0.525]{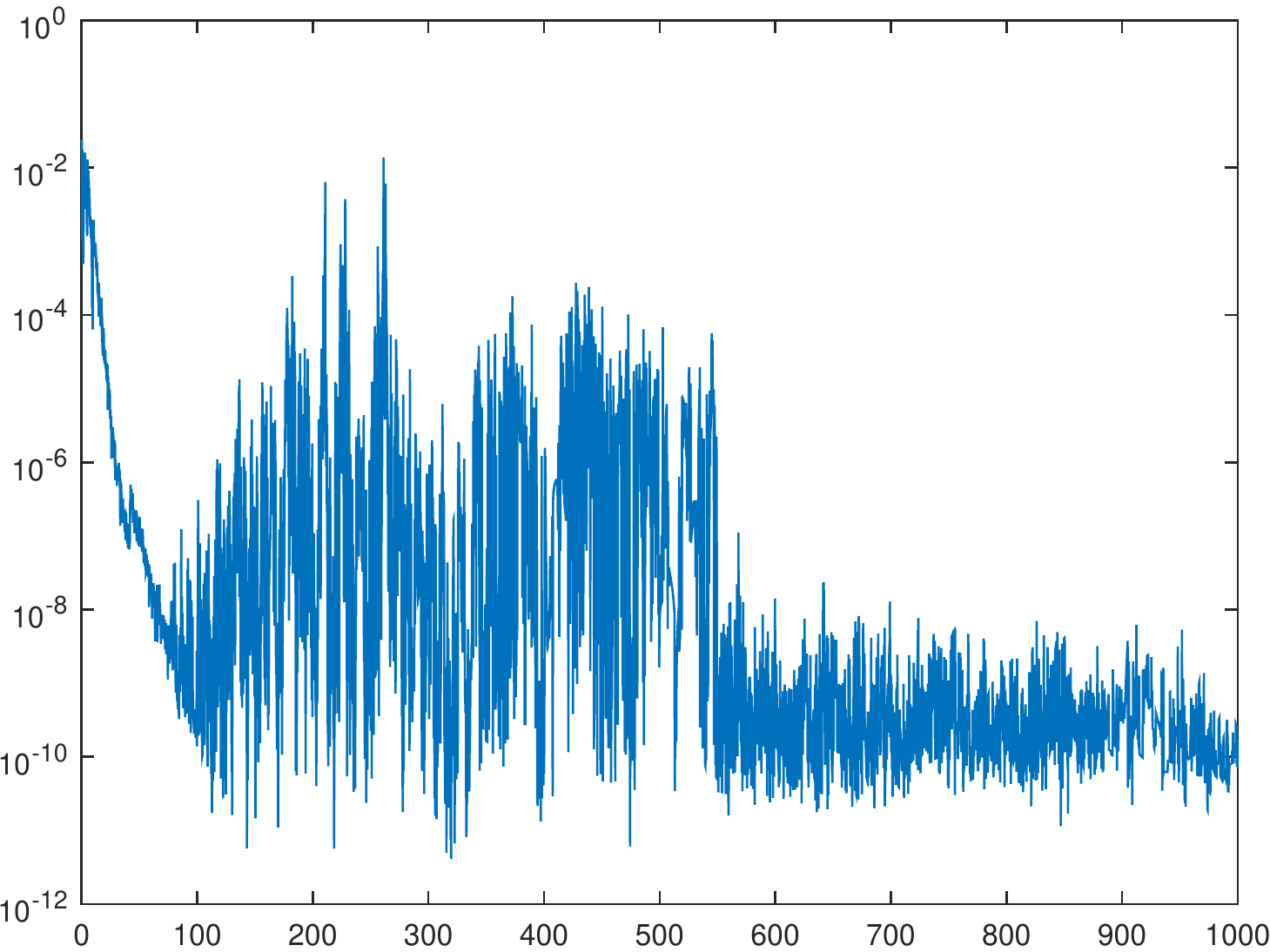}\\
		\includegraphics[scale=0.525]{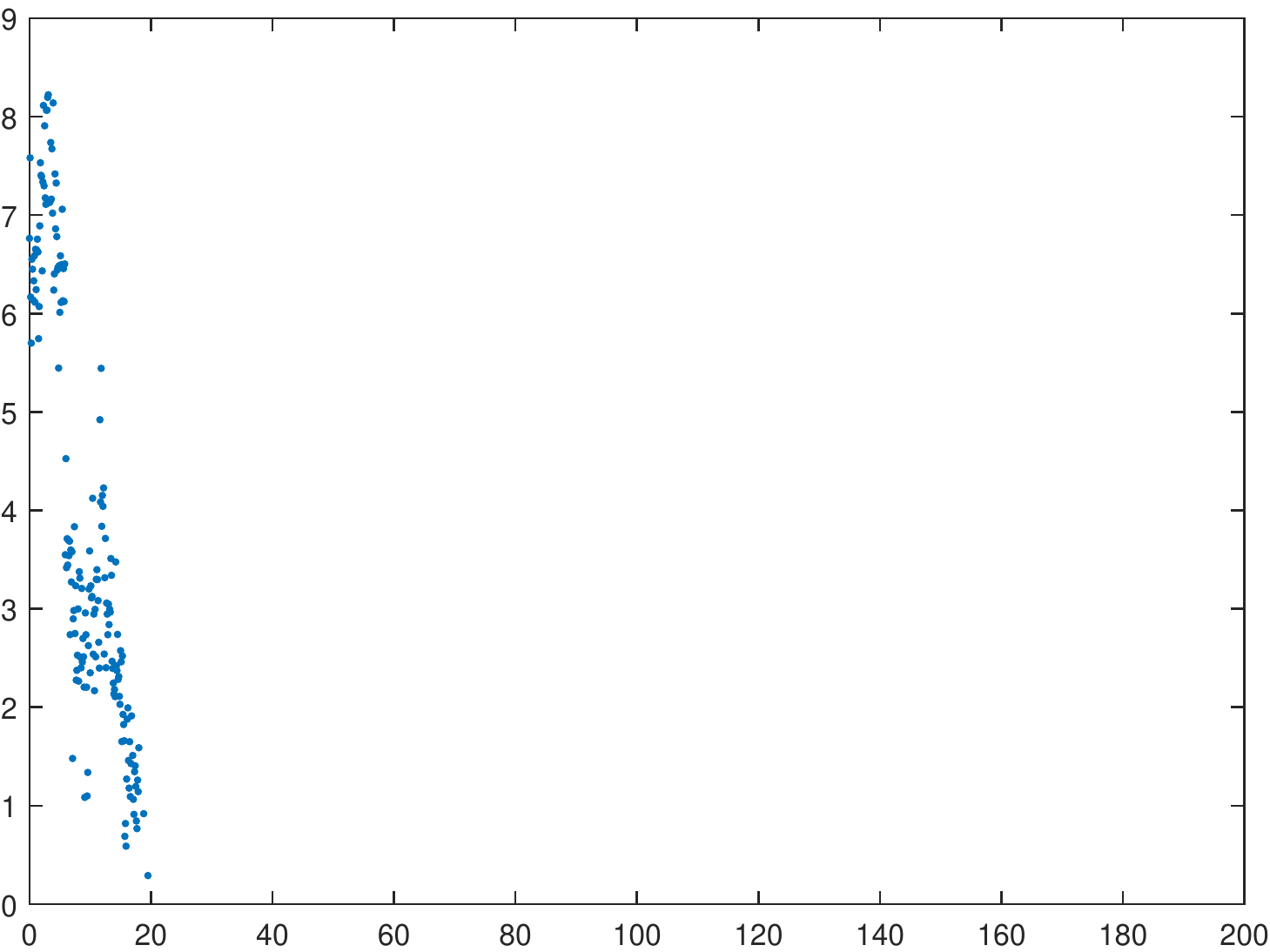}  & \includegraphics[scale=0.525]{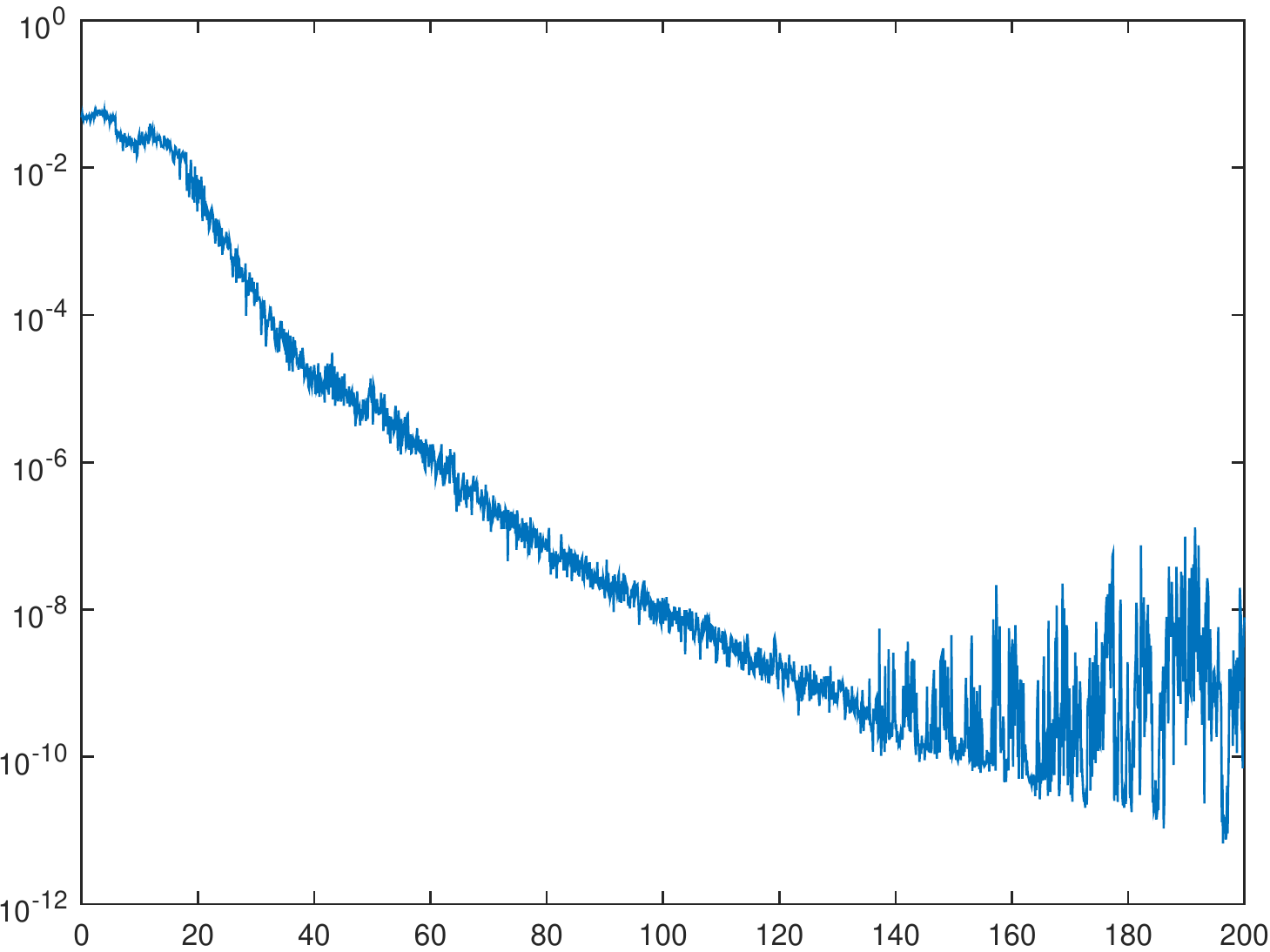}  \\
		\\
		\hline
	\end{tabular} 
	\caption{Plots of the largest imaginary part (left) and the relative error (right) in the computed eigenvalues of $H_{100}(a,b)$, horizontal axis varying values of $a$. Top row $b=50$, bottom row $b=100$.}
	\label{fig5}
\end{figure}

Next, we consider the general setting where we have two parameters $a$ and $b$, starting with the case where $n$ is even. Although the two parameters $a$ and $b$ occur symmetric in \eqref{eigevenab} and in the matrix \eqref{Hahneven} itself, there are some disparities in the numerical results. 
From the expression for the eigenvalues \eqref{eigevenab} we see that they are real when the parameters satisfy $a+b>-2$. 
However, 
\begin{itemize}
	\item when $a$ is a negative number less than $-2$, MATLAB finds eigenvalues with imaginary parts for almost all values of $b$.
	\item for negative values of $b$, no imaginary parts are found as long as $a+b>-2$.
	\item For positive values of $a$, eigenvalues with imaginary parts are found when $b$ gets sufficiently large, as illustrated in figure~\ref{fig4}. 
	\item For positive values of $b$ the opposite holds: eigenvalues with imaginary parts are found when $a$ is comparatively small, see figure~\ref{fig5}.
\end{itemize}



\begin{figure}[!htbp]
	\begin{tabular}{ccc}
		\hline \\
		\includegraphics[scale=0.525]{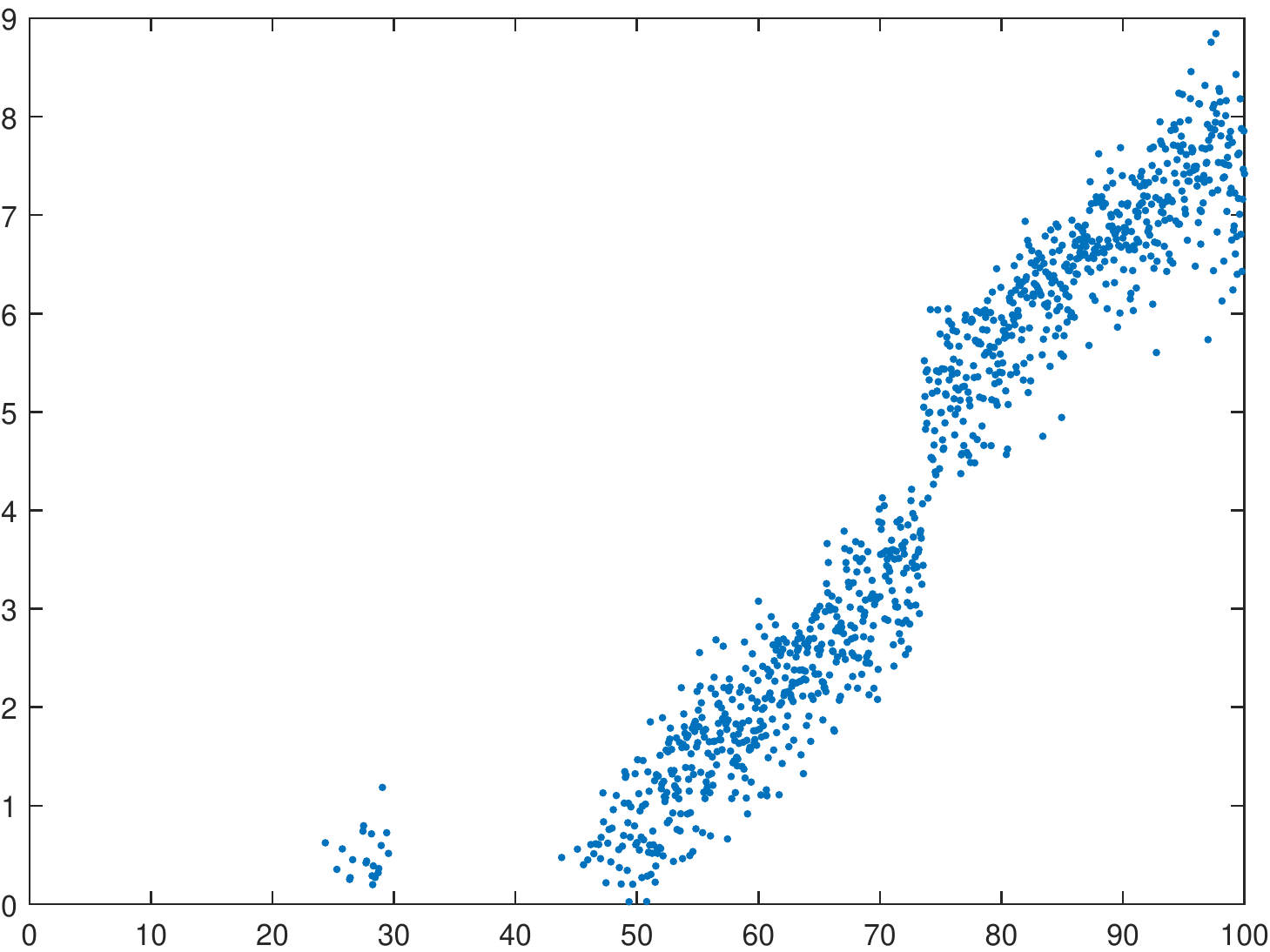} & \includegraphics[scale=0.525]{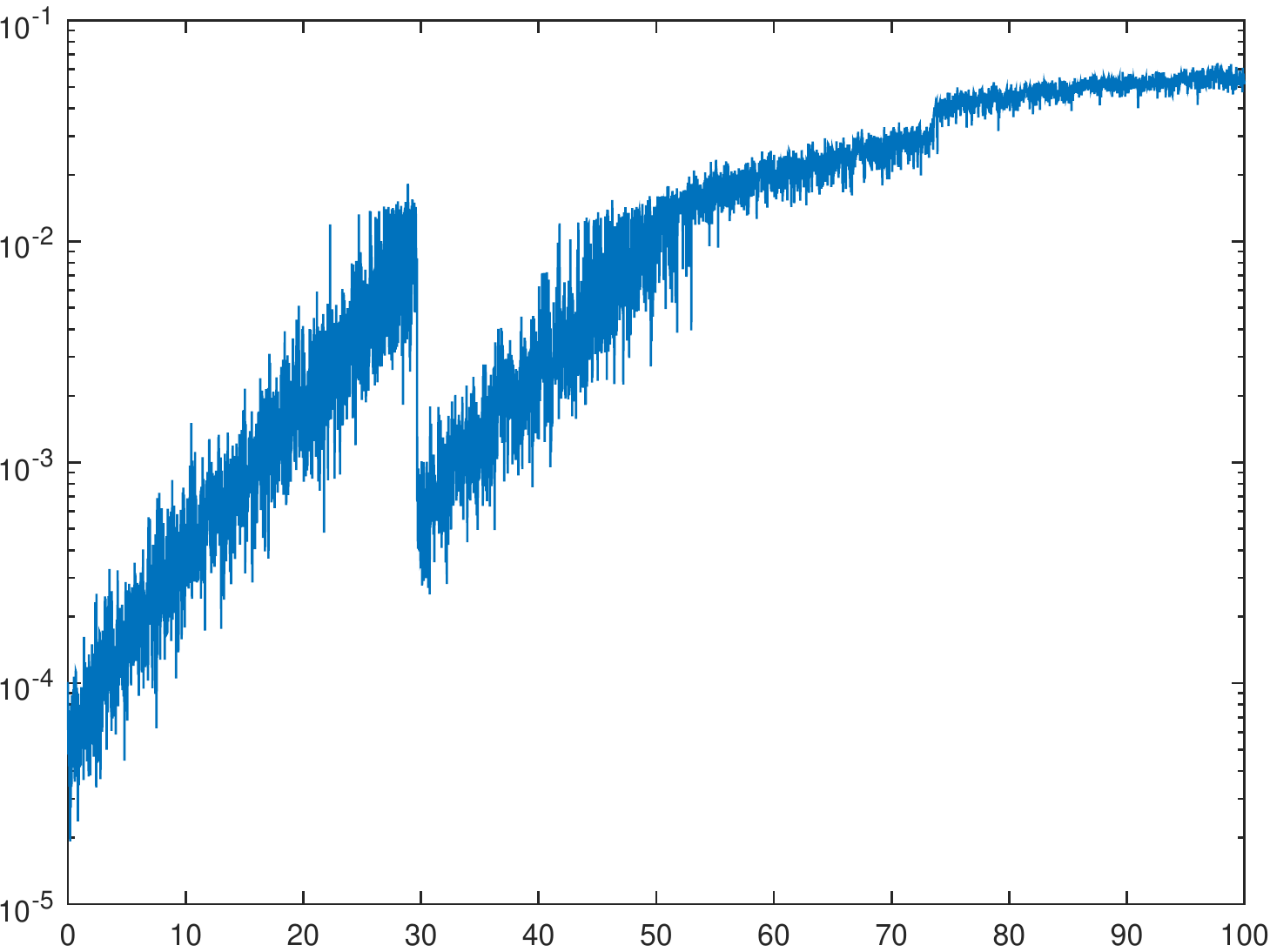} \\  \includegraphics[scale=0.525]{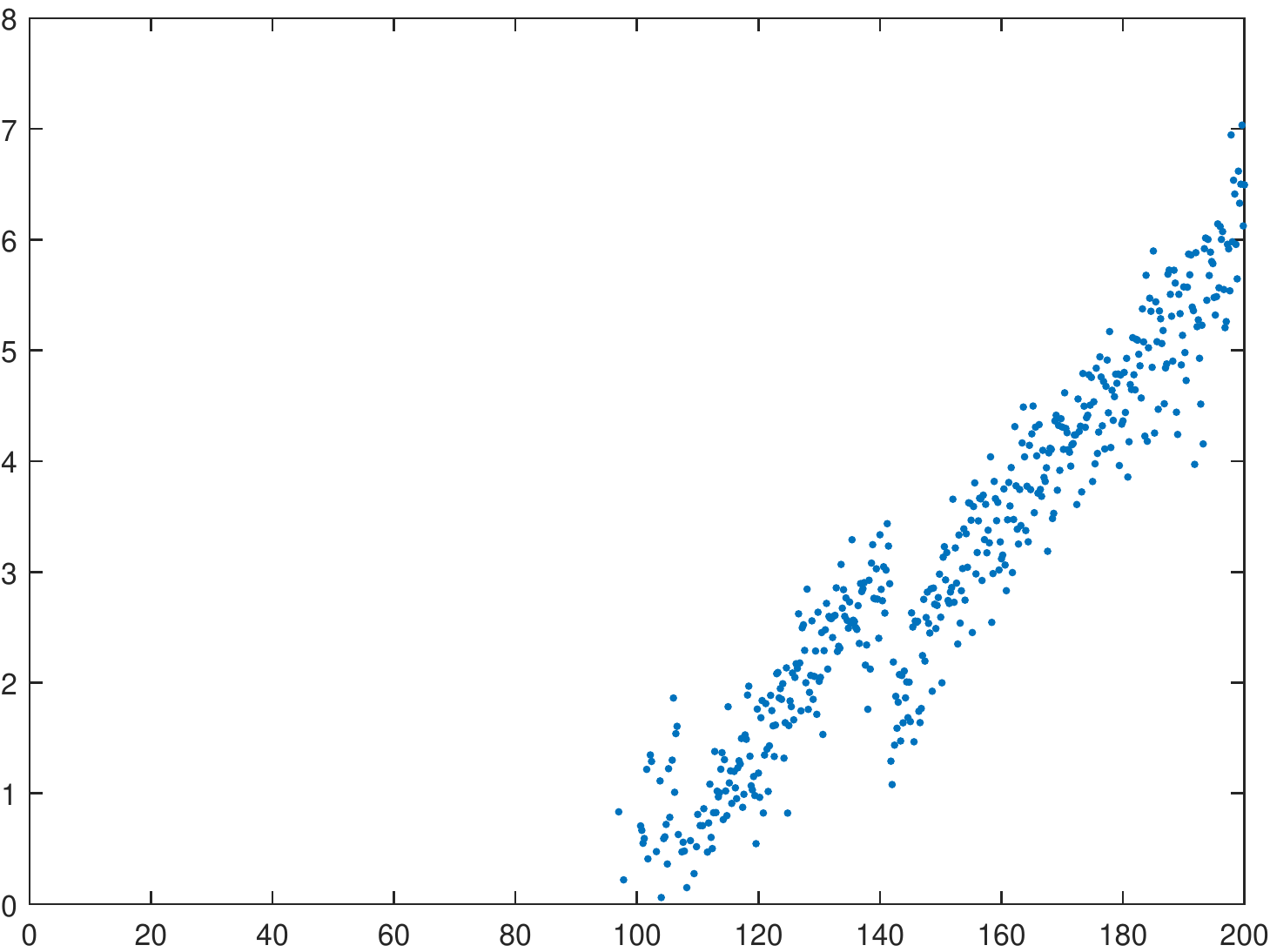} &
		\includegraphics[scale=0.525]{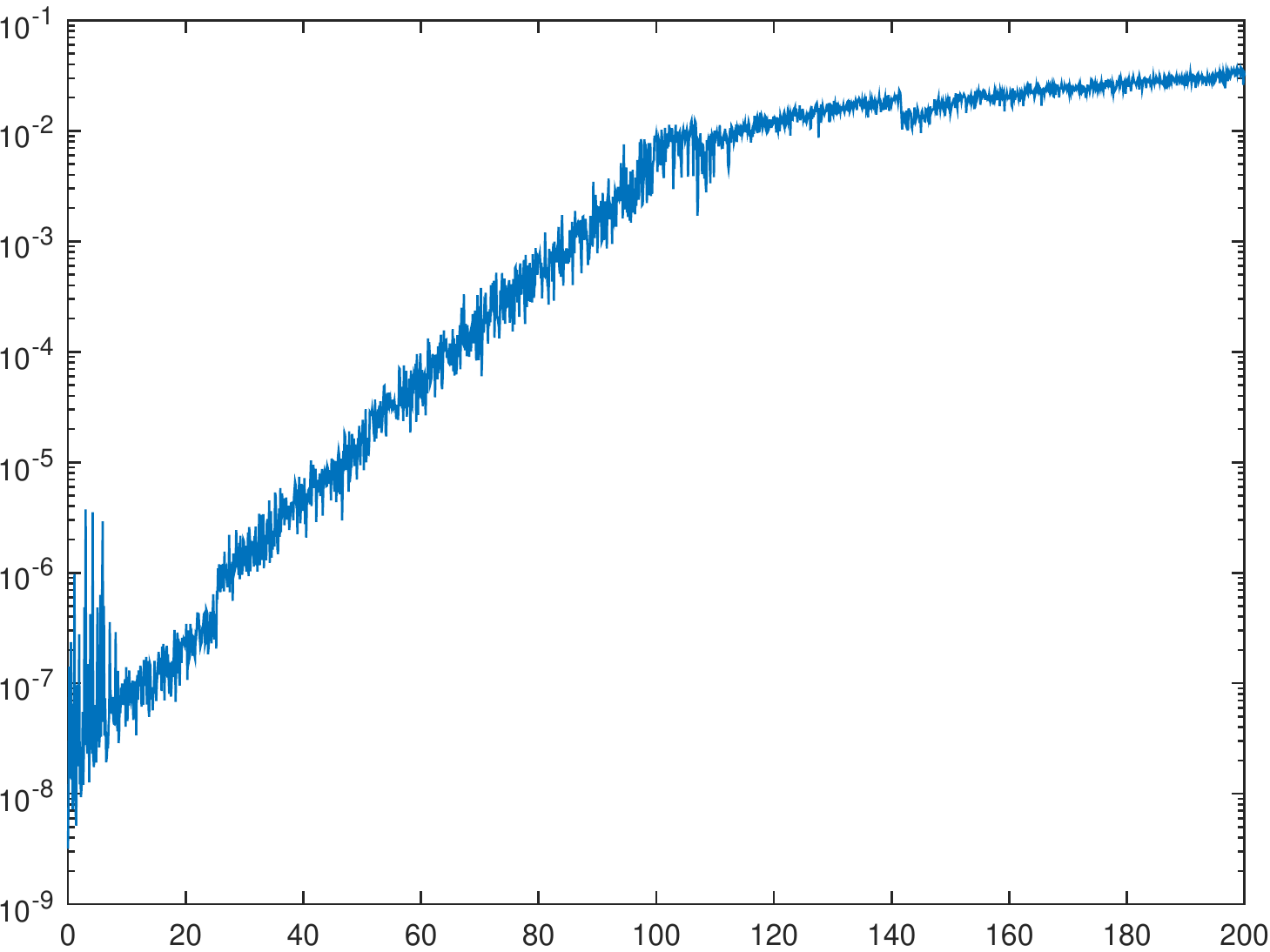} \\
		\\
		\hline
	\end{tabular} 
	\caption{Plots of the largest imaginary part (left) and the relative error (right) in the computed eigenvalues of $H_{101}(a,b)$, horizontal axis varying values of $b$. Top row $a=0$, bottom row $a=25$.}
	\label{fig6}
\end{figure}

In the case where $n$ is odd, similar results hold for positive parameter values for $a$ and $b$, as shown in figure~\ref{fig6} for example. For negative parameter values we have a different situation, as the eigenvalues \eqref{eigoddab} can become imaginary if the two factors have opposite sign. When $a<-n$ and $b<-n$, the eigenvalues \eqref{eigoddab} are real again and the behavior mimics that of the positive values of $a$ and $b$. The picture we get is a mirror image of figure~\ref{fig6}. 

The reason for this disparity between the seemingly symmetric parameters $a$ and $b$ is that the QR algorithm wants to get rid of subdiagonal entries in the process of creating an upper triangular matrix. As a consequence, the numerical computations are much more sensitive to large values of $b$ as it resides on the subdiagonal. 
This is showcased in figures~\ref{fig4}, \ref{fig5} and \ref{fig6}. 
Most important is the sensitivity on the extra parameters ($a$ or $a$ and $b$) which makes them appealing as test matrices.

It would be interesting future work if these new eigenvalue test matrices were to be used to test also numerical algorithms for computing eigenvalues designed specifically for matrices having multiple eigenvalues \cite{Galantai}, being tridiagonal \cite{Maeda}, or symmetric and tridiagonal \cite{Brockman}.

\newpage

\end{document}